\documentclass{amsart}
\newtheorem{Theorem}{Theorem}[section]
\newtheorem{Definition}[Theorem]{Definition}
\newtheorem{Proposition}[Theorem]{Proposition}

\newtheorem{Lemma}[Theorem]{Lemma}
\newtheorem{Corollary}[Theorem]{Corollary}
\theoremstyle{remark}
\newtheorem{Example}[Theorem]{Example}

\def\Xint#1{\mathchoice
{\XXint\displaystyle\textstyle{#1}}%
{\XXint\textstyle\scriptstyle{#1}}%
{\XXint\scriptstyle\scriptscriptstyle{#1}}%
{\XXint\scriptscriptstyle\scriptscriptstyle{#1}}%
\!\int}
\def\XXint#1#2#3{{\setbox0=\hbox{$#1{#2#3}{\int}$ }
\vcenter{\hbox{$#2#3$ }}\kern-.5\wd0}}

\def\dashint{\Xint-}
\def\il{\int\limits_}

\def\eps{\varepsilon}

\def\ovr{\overline}

\def\al{\alpha}

\def\th{\theta}

\def\dl{\delta}

\def\bd{\partial}
\def\lm{\lambda}

\def\sm{\setminus}
\def\sbs{\subset}
\def\sps{\supset}

\def\nea{\nearrow}
\def\sea{\searrow}
\def\wtl{\widetilde}

\def\supp{\operatorname{supp}}

\def\dist{\operatorname{dist}}

\def\be{\begin{enumerate}}
\def\ee{\end{enumerate}}
\def\bT{\begin{Theorem}}
\def\eT{\end{Theorem}}
\def\bP{\begin{Proposition}}
\def\eP{\end{Proposition}}
\def\bD{\begin{Definition}}
\def\eD{\end{Definition}}
\def\bE{\begin{Example}}
\def\eE{\end{Example}}
\def\bL{\begin{Lemma}}
\def\eL{\end{Lemma}}
\def\bC{\begin{Corollary}}
\def\eC{\end{Corollary}}
\def\A{{\mathcal A}}

\def\H{{\mathcal H}}

\def\F{{\mathcal F}}

\def\R{{\mathcal R}}
\def\rS{{\mathcal S}}
\def\T{{\mathcal T}}

\begin{document}
\title{Weak and strong limit values}
\author{Evgeny A. Poletsky}
\thanks{The author was supported by the NSF Grant DMS-0900877.}
\address{ Department of Mathematics,  215 Carnegie Hall,
Syracuse University,  Syracuse, NY 13244, eapolets@syr.edu}
\begin{abstract} The classical results about the boundary values of
holomorphic or harmonic functions on a domain $D$ state that under
additional integrability assumptions these functions have limits
along specific sets approaching boundary. The proofs of these
results are based on properties of smooth boundaries used to
define the approach regions and on estimates of representing
kernels along these regions.
\par This paper attempts to look at the situation when no
assumptions about the boundary smoothness are made and,
consequently, no natural definitions of approach regions could be
given.
\end{abstract}
\maketitle
\section{Introduction}
\par The theory of the boundary values of functions defined on a domain
$D\sbs\mathbb R^n$ and lying in some class $\rS$ tries to answer the
following questions:
\be\item Is there a trace operator $\T$ mapping $\rS$ into some class $\rS^b$ of
functions on $\bd D$ so that $\T u$ is neatly associated with $u$?
\item Is there a restoring operator $\R$ mapping $\rS^b\oplus\{\text{
    some information}\}$ so that $\T\R$ is an identity and $\R\T$ is
    the identity provided that ``some information" is supplied?
\ee
\par The most developed boundary values theory is for subclasses of (sub)harmonic functions. For bounded domains the classical Perron--Wiener--Brelot (PWB) method (see \cite{H}) provides harmonic solutions to the second problem for the class of functions integrable
with respect to the harmonic measure, while the Fatou--Na\"im--Doob
Theorem (see \cite{AG}) asserts the existence of minimal fine limits of non-negative harmonic functions almost everywhere on the boundary.
However, no connections between the latter limits and the given
function were established for general domains.
\par When $n=2$ and $D$ is the unit disk the work of P. Fatou and F.
Riesz completely answered these questions when $\rS$ is a Hardy space either of harmonic functions  or holomorphic functions. It was shown that for a function $h$ in such a space non-tangential limits exist almost everywhere at the boundary with respect to the surface length and play the role of the trace $\T h$. The Poisson integral plays the role of $\R$.  J. E. Littlewood expanded this work to the class of subharmonic functions with harmonic majorants using the Laplacian as ``some information" and the Riesz decomposition formula as $\R$.
\par In 1936 I. I. Privalov and P. K. Kuznetsov lifted Littlewood's
results to $n=3$ and $D$ be a ball. E. D. Solomentsev in \cite{S}
generalized their work to any $n$ and $D$ be a $C^2$-domain. In
\cite{Da} B. E. J. Dahlberg expanded this result to $C^{1,1}$-domains.
D. S. Jerison and C. E. Kenig continued this work in \cite{JK1} and for a larger class of $L^p_1$-domains showed that under an integrability
condition on the boundary data $\phi$ the PWB solution is given by the
Poisson formula and has non-tangential limits equal to $\phi$ almost
everywhere on the boundary. They also proved that bounded harmonic
functions on such domains have non-tangential limits almost everywhere
at the boundary with respect to the surface area.
\par In \cite{JK2} D. S. Jerison and C. E. Kenig generalized their results to non-tangentially accessible domains. However, the two major changes of the classical theory were made in this paper. First of all, the surface area was abandoned and replaced by the harmonic measure $\mu$ and, secondly, the Hardy spaces were defined as sets of functions for which the maximal non-tangential function belongs to $L^p(\mu)$. This approach was mostly used in recent publications. 
\par Similar results for subharmonic functions seem to be non-existent. Moreover, in \cite{Da} B. E. J. Dahlberg showed that there are a $C^{1,\al}$-domain $D\sbs\mathbb R^2$ ($0<\al<1$) and a negative subharmonic function on $D$ which has no limits along normals on a set of positive length. So the validity of the the Riesz decomposition formula was not clear.
\par The situation was even worse for plurisubharmonic functions, where the Lelong--Jensen formula (see \cite{D}) replaces the Riesz decomposition formula. Roughly speaking, to get good boundary values results one, firstly, establishes them on compact subdomains and then exhausts $D$ by open sets $D_j\sbs\sbs D$, $j=1,2,\dots$. If the boundary is nice the exhaustion can be made even nicer. The formulas representing functions at $x\in D$ (e.g., Riesz or Lelong--Jensen formulas) involve an integral of the function over the boundary
$S_j=\bd D_j$ over some measure $\mu_{xj}$ and an integral over $D_j$.
\par While it is easy to establish that the volume integrals converge
to similar integrals as $j\to\infty$, the surface integrals present
the major difficulty. In the classical potential theory the measures  $\mu_{xj}$ are the surface areas times the Poisson kernel. The behavior of the Poisson kernel on smooth domains is well studied and it allows to prove that the surface integrals converge to similar integrals of pre-existing boundary values.
\par However, an attempt to replicate this for subharmonic functions on
non-smooth domains or plurisubharmonic functions on hyperconvex domains
fails due to the absence of knowledge about the measures $\mu_{xj}$. In the case of plurisubharmonic functions except of \cite{D} only \cite{BPT} addresses the latter problem for strongly convex domains with smooth
boundary.
\par In this paper we suggest an approach which, firstly, abandons the non-tangential limits and, secondly, returns to the classical definition of Hardy spaces. To understand the situation we start with the general problem of different type convergence for sequences $\{\phi_j\mu_j\}$, where $\mu_j$ is a measure on $S_j$. Standard integrability conditions show that the sequence $\{\phi_j\mu_j\}$ has a subsequence converging weak-$*$ to a measure $\phi_*\mu_x$ and we restrict our attention to weak-$*$ converging sequences. The functions $\phi_*$ can be considered as the weak limit values. The main disadvantage of these functions $\phi_*$ is their bad correlation with products. It is not true, in general, that $(\phi\psi)_*=\phi_*\psi_*$ - the identity one needs to prove the integral formulas (see examples in Section \ref{S:wlv}).
\par However, we establish that if the sequence $\{|\phi_j|^p\mu_j\}$ also has the weak-$*$ limit $\nu$, then $\nu\ge|\phi_*|^p\mu$. Surprisingly,
$\nu=|\phi_*|^p\mu$ if and only if $\phi_*$ has much stronger
properties and we call it the strong limit values and denote by
$\phi^*$ (for the precise definition see Section \ref{S:slv}. This
theory is developed in Sections \ref{S:wlv} and \ref{S:slv}. In
particular, we show that if $\psi_*$ exists then
$\phi^*\psi_*=(\phi\psi)_*$.
\par In Section \ref{S:ssii} we address the question when the
strong limit values exist or when $\nu=|\phi_*|^p\mu$. We consider a
space of functions with the weak limit values and satisfying an
integral inequality (\ref{e:tc}). For example, if the space in question
is the space of subharmonic functions, then the inequality is the
classical estimate of the value of a function $\phi$ at some point
through the convolution of $\phi_*$ and the Poisson kernel.
\par Under mild conditions on the kernel we prove that the functions
in such spaces have the strong limit values. The section also contains
several results showing when these mild conditions hold.
\par In Section \ref{S:bv} we look at functions defined on $D$ and define boundary values as the strong limit values for all possible sequences
of exhaustions. We give some sufficient conditions for functions to
have boundary values, establish a general result allowing to expand
integral representations and show that functions which have classical
radial limits also have boundary values in our sense.
\par Section \ref{S:bvhf} contains applications of exposed theory to
harmonic and subharmonic functions on a general regular domain $D$. We introduce the Hardy spaces $\rS^p(D)$ of (sub)harmonic functions on $D$ and show that any function in $L^p(\bd D,\mu_{x_0})$ is the trace of a function in $\rS^p(D)$ and the Poisson integral serves as restoring operator. We also prove that the Riesz decomposition formula is valid for a function $u\in\rS^p(D)$ if and only if $u$ has the boundary values. 
\par The last section \ref{S:sopf} contains the expansion of
Lelong--Jensen formula to $L^p$-classes of plurisubharmonic functions
on strongly pseudoconvex domains and the polydisk. This case was also
studied in \cite{CK}.
\par The author is grateful to Michael Stessin for his strong
encouragement to write this paper and to the referee whose suggestions improved the exposition.
\section{Weak limit values}\label{S:wlv}
\par Let $K$ be a compact metric space, and let $M=\{\mu_j\}$
be a sequence of regular Borel measures on $K$ converging
weak-$*$ in $C^*(K)$ to a finite measure $\mu$. We denote the
set $\supp\mu_j$ by $K_j$ and $\supp\mu$ by $K_0$. Let
$\phi=\{\phi_j\}$ be a sequence of Borel functions $\phi_j$ on
$K_j$. We let
\[\|\phi\|_{L^p(M)}=\limsup_{j\to\infty}\|\phi_j\|_{L^p(K_j,\mu_j)}.\]
\par In general, the weak-$*$ limit of measures $\phi_j\mu_j$ need
not to be absolutely continuous with respect to $\mu$. But as the
following lemma shows this is the case when $\|\phi\|_{L^p(M)}<\infty$
for some $p>1$.
\bL\label{L:bv1} If $\|\phi\|_{L^p(M)}\le A$, $p>1$, and the measures $\phi_j\mu_j$
converge weak-$*$ to a measure $\mu'$ on $K$, then there is a
function $\phi_*\in L^p(K_0,\mu)$ such that $\mu'=\phi_*\mu$
and $\|\phi_*\|_{L^p(K_0,\mu)}\le A$.
\eL
\begin{proof} First of all, we show that $\supp\mu'=\supp\mu$, i.e., for
any $h\in C(K)$
\[\mu'(h)=\il Kh\,d\mu'=\il{K_0}h\,d\mu'.\]
For this for $\dl>0$ we take the closed $\dl$-neighborhood $V$ and
the open $2\dl$-neighborhood $W$ of $K_0$. Let $f_1$ be a
non-negative continuous function on $K$ taking values between 0
and 1, which is equal to 1 on $V$ and whose support lies in $W$.
Let $f_2=1-f_1$. Then
\[\lim_{j\to\infty}\il {K_j}f_2h\,d\mu_j=\il {K_0}f_2h\,d\mu=0.\]
By H\"older's inequality
\[\left|\il {K_j}f_2h\phi_j\,d\mu_j\right|
\le\|\phi_j\|_{L^p(K_j,\mu_j)}\left(\il
{K_j}f_2h\,d\mu_j\right)^{1/q},\] where $1/p+1/q=1$. Thus
\[\mu'(h)=\mu'(f_2h)+\mu'(f_1h)=\lim_{j\to\infty}\il
Kf_2h\phi_j\,d\mu_j+\mu'(f_1h)=\mu'(f_1h).\]
\par Hence for every open neighborhood $Y$ of $K_0$
\[\mu'(h)=\il Yh\,d\mu'\]
and this implies that
\[\mu'(h)=\il Kh\,d\mu'=\il{K_0}h\,d\mu'.\]
\par By H\"older's inequality
\[\left|\il{K_j}h\phi_j\,d\mu_j\right|
\le \|\phi_j\|_{L^p(K_j,\mu_j)}\|h\|_{L^q(K_j,\mu_j)}.\] Hence
\[\left|\il{K}h\,d\mu'\right|=
\lim_{j\to\infty}\left|\il{K_j}h\phi_j\,d\mu_j\right|\le
A\|h\|_{L^q(K_0,\mu)}.\] So the functional
\[\mu'(h)=\il{K}h\,d\mu'=\il{K_0}h\,d\mu'\] on $C(K)$ admits an unique
extension to a continuous functional on $L^q(K_0,\mu)$. Thus there
is a function $\phi_*\in L^p(K_0,\mu)$ such that
$\|\phi_*\|_{L^p(K_0,\mu)}\le A$ and
\[\il{K_0}h\,d\mu'=\il {K_0}h\phi_*\,d\mu\] for any $h\in C(K)$. Therefore,
$\mu'=\phi_*\mu$.
\end{proof}
\par If the measures $\{\phi_j\mu_j\}$ converge weak-$*$ to a measure
$\phi_*\mu$, then the function $\phi_*$ will be called the {\it weak
limit values } of $\phi$. We will denote by $\A(M)$ the space of all
sequences $\phi$ of Borel functions $\phi_j$ on $K_j$  which have weak
limit values and by $\A^p(M)$ those sequences $\phi$ in $\A(M)$ for
$\|\phi\|_{L^p(M)}<\infty$.
\par While the weak-$*$ convergence of measures $\{\phi_j\mu_j\}$ frequently
occurs it does not correlate well with algebraic operations. First of
all, as the following example shows, it does not imply the weak-$*$
convergence of the sequence $\phi^p=\{|\phi_j|^p\mu_j\}$.
\par {\bf Example:} Let $K=K_0=[0,1]$ and all measures $\mu_j$ be equal to
the Lebesgue measure $\mu$ on $[0,1]$. For each $j$ we split $[0,1]$
into intervals $I_{jk}=\{x:k/j\le x<(k-1)/j\}$, $k=0,\dots,j-1$. When
$j$ is even we let $\phi_j(x)=0$ if $x\in I_{jk}$ and $k$ is even and
we let $\phi_j(x)=2$ if $x\in I_{jk}$ and $k$ is odd. When $j$ is odd
then we let $\phi_j(x)=1$ for all $x\in[0,1]$. Then the measures
$\phi_j\mu_j$ converge weak-$*$ to $\mu$, while the measures
$\phi^2_j\mu_j$ do not converge weak-$*$. For even $j$ they converge to
$2\mu$ and for odd ones to $\mu$.
\par Even if the sequence $\phi^p=\{|\phi_j|^p\mu_j\}$
weak limit values it is not true that $|\phi_*|^p=(|\phi^p|)_*$ as the
example below shows.
\par {\bf Example:} We take measures and intervals from the previous example
and for all $j$ we let $\phi_j(x)=a>0$ if $x\in I_{jk}$ and $k$ is even
and we let $\phi_j(x)=b>0$ if $x\in I_{jk}$ and $k$ is odd. Then the
measures $\phi_j\mu_j$ converge weak-$*$ to $\phi_*\mu$, where
$\phi_*\equiv(a+b)/2$, while the measures $\phi^2_j\mu_j$ converge
weak-$*$ to $\nu=(a^2+b^2)\mu/2$.
\par Note that $\phi_*^2\mu\le\nu$ and the theorem below shows that this
crucial observation is always true.
\bT\label{T:co} If $\phi\in\A^p(M)$, $p>1$, and the measures
$|\phi_j|^p\mu_j$ converge weak-$*$ to a measure $\nu$ on $K$,
then $\nu\ge|\phi_*|^p\mu$.
\eT
\begin{proof} Let $A=\|\phi_*\|_{L^p(K,\mu)}$ and let
$G(a,b)=\{x\in K_0:\,a\le\phi_*(x)<b\}$. Suppose that $b>a>0$ and let
$X\sbs G(a,b)$ be a Borel set. Suppose that $\mu(X)=m>0$. We fix
$\eps>0$, such that $a\mu(X)-A\eps^{1/q}>0$, $1/p+1/q=1$, and take an
open set $O$ and a closed set $C$ such that $C\sbs X\sbs O$,
$a\mu(C)-A\eps^{1/q}>0$, $\nu(O\sm C)<\eps$ and $\mu(O\sm C)<\eps$. Let
$f$ be a continuous function on $K$ equal to 0 on $K\sm O$, equal to 1
on $C$ and taking values between 0 and 1 elsewhere.
\par Then
\begin{equation}\begin{aligned}
\nu(O)&\ge\il Kf\,d\nu=\lim_{j\to\infty}\il Kf|\phi_j|^p\,d\mu_j
\ge\lim_{j\to\infty}\left|\il Kf\phi_j\,d\mu_j\right|^p\left(\il
Kf\,d\mu_j\right)^{1-p}\notag\\&=\left|\il
Kf\phi_*\,d\mu\right|^p\left(\il
Kf\,d\mu\right)^{1-p}.\notag\end{aligned}\end{equation} Now
\[\il Kf\phi_*\,d\mu=\il Cf\phi_*\,d\mu+\il{O\sm C}f\phi_*\,d\mu\ge
a\mu(C)+\il{O\sm C}f\phi_*\,d\mu\] and
\[\left|\il{O\sm
C}f\phi_*\,d\mu\right|\le\|\phi_*\|_{L^p(K,\mu)}\|f\|_{L^q(K,\mu)}\le A\eps^{1/q}.\]
Hence
\[\left|\il
Kf\phi_*\,d\mu\right|^p\ge\left(a\mu(C)-A\eps^{1/q}\right)^p.\]
Thus
\[\nu(O)\ge\left(a\mu(C)-A\eps^{1/q}\right)^p\mu^{1-p}(C).\]
Letting $\eps$ go to 0 we see that $\nu(X)\ge a^p\mu(X)$. Hence
\[\il X\phi_*^p\,d\mu\le b^p\mu(X)\le\frac{b^p}{a^p}\nu(X).\]
\par If $b<a<0$ then replacing $\phi$ by $-\phi$ we see that
\[\il X|\phi_*|^p\,d\mu\le |b|^p\mu(X)\le\left|\frac{b}{a}\right|^p\nu(X).\]
Note that if $\mu(X)=0$ both inequalities above are trivially
true.
\par If $X$ is any Borel set in $K_0\cap\{\phi_*\ge0\}$, we
denote by $X_\eps$, $\eps>0$, the intersection of $X$ and the
set $\{\phi_*>\eps\}$, where $\eps>0$. Let us take some
$\dl>0$, form a sequence $a_0=\eps$, $a_{k+1}=(1+\dl)^{1/p}a_k$
and let $X_\eps^k=X_\eps\cap G(a_k,a_{k+1})$. Then
\[\il {X_\eps}\phi_*^p\,d\mu=\sum_{k=0}^\infty\il
{X^k_\eps}\phi_*^p\,d\mu\le(1+\dl)\nu(X_\eps).\] Thus
\[\il {X_\eps}\phi_*^p\,d\mu\le\nu(X).\]
\par If $X'_\eps=X\sm X_\eps$ then
\[\il {X'_\eps}\phi_*^p\,d\mu\le\eps^p\mu(X).\] Thus
\[\il {X}\phi_*^p\,d\mu\le\nu(X)+\eps^p\mu(X)\] and it shows that
\[\il {X}\phi_*^p\,d\mu\le\nu(X).\]
\par If $X\sbs\supp\mu\cap\{\phi_*<0\}$ then a consideration of the
function $-\phi$ shows that
\[\il {X}|\phi_*|^p\,d\mu\le\nu(X).\] Thus $\nu\ge|\phi_*|^p\mu$.
\end{proof}
\par We finish this section with a couple of technical lemmas.
\bL\label{L:rcs} Let $\{\mu_j\}$ be a sequence of regular Borel measures on $K$
converging weak-$*$ to $\mu$ and let $A$ be a closed set in $K$. Then
for every $\eps>0$ there is an open set $O$ containing $A$ and $j_0$
such that $\mu_j(O)<\mu(A)+\eps$ when $j\ge j_0$.
\eL
\begin{proof} Let us fix some $\dl>0$ and $0<a<1$ whose precise values
will be determined later and take a continuous function $f$ on $K$
which is equal to 1 on $A$, $0\le f<1$ elsewhere and
$\mu(f)<\mu(A)+\dl$. Then there is $j_0$ such that
$\mu_j(f)<\mu(f)+\dl$ when $j\ge j_0$. Since $\mu_j(f)\ge
a\mu_j(\{f>a\})$ we see that
\[\mu_j(\{f>a\})<\frac{\mu(f)+\dl}a<\frac{\mu(A)+2\dl}a.\]
If we take $a$ and $\dl$ satisfying the inequality
$2\dl<a\eps+(a-1)\mu(A)$ and let $O=\{f>a\}$ then
$\mu_j(O)<\mu(A)+\eps$.
\end{proof}
\par This lemma, in general, does not hold for non-closed sets.
For example, let $\mu$ be the Lebesgue measure on $[0,1]$ and
\[\mu_j=\frac1n\sum_{k=0}^n\dl_{k/j}.\]
Then $\mu_j$ converge weak-$*$ to $\mu$, $\mu(\mathbb Q)=0$ but
$\mu_j(\mathbb Q)\ge1$.
\par To resolve this problem we introduce the {\it regular}
weak-$*$ convergence. Namely, we say that measures $\{\mu_j\}$ converge
weak-$*$ to $\mu$ {\it regularly } if they converge weak-$*$ and for
any $\mu$-measurable set $A$ in $K_0=\supp\mu$ and every $\eps>0$ there
is an open set $O$ containing $A$ and $j_0$ such that
$\mu_j(O)<\mu(A)+\eps$ when $j\ge j_0$.
\par Usually, it is easy to verify the regularity of a weak-$*$
convergence. If not then the lemma below gives a criterion which
is especially useful for boundary values.
\bL\label{L:rc} Suppose that in assumptions of Lemma \ref{L:rcs}
$K_j\cap K_0=\emptyset$ for all $j$. Then the measures
$\{\mu_j\}$ converge weak-$*$ to $\mu$ regularly.
\eL
\begin{proof} We already know that the lemma holds for all closed set.
Let $A\sbs K_0$ be a relatively open set. We look at the distance
function $d(x)=\dist(x,\bd A)$ on $K_0$ and note that
$\mu(\{d(x)=a\})>0$ only for countably many values of $a$.
Choosing a sequence of $a_k\sea 0$ such that $\mu(\{d(x)=a_k\})=0$
we let $A_k=\{x\in A:\,a_k\le d(x)\le a_{k+1}\}$. Then the sets
$A_k$ are closed, $A=\cup A_k$ and $\mu(A)=\sum\mu(A_k)$.
\par For $\eps>0$ and any $k$ we find $j_k$ and an open set $O_k$ such
that $A_k\sbs O_k$ and $\mu_j(O_k)<\mu(A_k)+2^{-k}\eps$ when
$j>j_k$. Let $O'_k=O_k\sm\cup_{j=1}^{j_k}K_j$. Then the sets
$O'_k$ are open, contain $A_k$ and
$\mu_j(O'_k)<\mu(A_k)+2^{-k}\eps$ for all $j$. If $O=\cup O'_k$
then $\mu_j(O)<\mu(A)+\eps$ for all $j$.
\par If $A$ is any $\mu$-measurable set in $\supp\mu$ then we find
a relatively open set $A'\sbs K_0$ such that
$\mu(A')<\mu(A)+\eps/2$ and then an open set $O$ in $K$ containing
$A'$ such that $\mu_j(O)<\mu(A')+\eps/2$ for all $j$. Then
$\mu_j(O)<\mu(A)+\eps$.
\end{proof}
\section{Strong limit values}\label{S:slv}
\par We will need a more precise H\"older's inequality. For $x\ge0$ and $p>1$ let
\begin{equation}\label{e:phi}
\Phi_p(x)=x^p-p(x-1)-1.
\end{equation}
Clearly, $\Phi_p(x)\ge0$ and $\Phi_p(x)=0$
if and only if $x=1$.
\par If $p\ge2$ then $\Phi_p(x)\ge|x-1|^p$. Indeed, if $x>1$ then
differentiating both sides we come to an evident inequality
\[x^{p-1}-(x-1)^{p-1}\ge 1.\] If $0\le x\le1$, then the same trick
leads to another evident inequality
\[x^{p-1}+(1-x)^{p-1}\le 1.\]
\par Let
\[\bar f=\dashint_Kf\,d\mu=\frac1{\mu(K)}\il Kf\,d\mu.\]
Replacing $x$ with $f/\bar f$ in (\ref{e:phi}) and integrating
both sides we get the following proposition.
\bP\label{P:phi} Let $(K,\mu)$ be a measure space, $0<\bar\mu=\mu(K)<\infty$,
and let $f$ be a non-negative measurable function on $K$ with
$\|f\|_{L^1(K,\mu)}<\infty$. If $p>1$ then
\begin{equation}\label{e:hi}\left(\dashint_Kf\,d\mu\right)^{-p}
\dashint_Kf^p\,d\mu=1+
\dashint_K\Phi_p\left(\frac{f}{\bar f}\right)\,d\mu.\end{equation}
\par If $p\ge2$ then
\[\left(\dashint_Kf\,d\mu\right)^{-p}\dashint_Kf^p\,d\mu\ge1+
\dashint_K|f/\bar f-1|^p\,d\mu.\]
\eP
\par We will need the following immediate consequence of this
proposition.
\bL\label{L:chi} In assumptions of Proposition \ref{P:phi} for $p>1$ and
$c>0$ there is a constant $\al(p,c)$ such that if the left side of
(\ref{e:hi}) is smaller than $1+\eps$, $\eps>0$, then $\mu(\{|f/\bar
f-1|>1+c\})<\al(p,c)\eps\mu(K)$.
\eL
\par We say that a sequence {\it $\phi\in\A(M)$ has the strong limit values on
$K_0$ with respect to $M$} if there is a $\mu$-measurable function
$\phi^*$ on $K_0$ such that for any $b>a$ and any $\eps,\dl>0$ there is
$j_0$ and an open set $O\sbs K$ containing $G(a,b)=\{x\in
K_0:\,a\le\phi^*(x)<b\}$ such that
\begin{equation}\label{e:bv}\mu_j(\{\phi_j<a-\eps\}\cap
O)+\mu_j(\{\phi_j>b+\eps\}\cap O)<\dl\end{equation} when $j\ge
j_0$. The function $\phi^*$ will be called the {\it strong
limit values} of $\phi$.
\par Let us indicate some properties of strong limit values.
\bT\label{T:pslf} Suppose that $\phi$ has the strong limit values
on $K_0$ equal to $\phi^*$. Then:
\be\item any two choices of $\phi^*$ coincide $\mu$-a.e., the sequences $c\phi$ and
$|\phi|^p$ have strong limit values and $(c\phi)^*=c\phi^*$ and
$(|\phi|^p)^*=|\phi^*|^p$;
\item if the sequence $\{\mu_j\}$ converges weak-$*$ regularly
    to $\mu$ and a sequence  $\psi\in\A(M)$ has the strong limit values
    $\psi^*$, then  the sequences $\phi+\psi$ and $\phi\psi$ have the strong
    limit values and $(\phi+\psi)^*=\phi^*+\psi^*$ and
    $(\phi\psi)^*=\phi^*\psi^*$.
\ee
\eT
\begin{proof} 1) Suppose that $\phi$ has two different strong limit values
$\phi^*_1$ and $\phi^*_2$. Suppose that for some $\eps>0$ there are
$a<b$ and $c<d$ such that $c>b+2\eps$ and
$\mu(\{a\le\phi_1^*<b\}\cap\{c\le\phi_2^*<d\})=\al>0$. Let us take open
sets $O_1$ containing $\{a\le\phi_1^*<b\}$ and $O_2$ containing
$\{c\le\phi_2^*<d\}$ such that
\[\mu_j(\{\phi_j<a-\eps\}\cap O_1)+\mu_j(\{\phi_j>b+\eps\}\cap O_1)<\al/4\]
and
\[\mu_j(\{\phi_j<c-\eps\}\cap O_2)+\mu_j(\{\phi_j>d+\eps\}\cap O_2)<\al/4\]
for large $j$.
\par If $O=O_1\cap O_2$ then $\mu_j(O)>\al/2$ for large $j$. Hence
\[\mu_j(\{c-\eps\le\phi_j\le d+\eps\}\cap O)>\al/4\] and,
consequently, $\mu_j(\{b+\eps<\phi_j\}\cap O)>\al/4$ and we get a
contradiction. Hence,
$\mu(\{a\le\phi_1^*<b\}\cap\{c\le\phi_2^*<d\})=0$ when
$c>b+2\eps$.
\par Since $\eps>0$ is arbitrary we see that
$\mu(\{a\le\phi_1^*<b\}\cap\{c\le\phi_2^*<d\})=0$ when $c>b$. Thus
$\phi_2^*\le b$ $\mu$-a.e. on $\{a\le\phi_1^*<b\}$ and this shows that
$\phi^*_2\le\phi_1^*$ $\mu$-a.e. By symmetry, $\phi^*_2\le\phi_1^*$
$\mu$-a.e.
\par The facts that $(c\phi)^*=c\phi^*$ and $(|\phi|^p)^*=|\phi^*|^p$ are trivial.
\par 2) We will prove only that $(\phi+\psi)^*=\phi^*+\psi^*$.
Other claims can be proved similarly.
\par For given numbers $\eps,\dl>0$ we find $c<d$, such that
$\mu(X)<\dl/2$, where
\[X=\{\phi^*<c\}\cup\{\psi^*<c\}\cup\{\phi^*>d\}\cup\{\psi^*>d\}.\]
Then we split the interval $[c,d]$ into consecutive intervals
$[a_k,a_{k+1}]$ such that $a_{k+1}-a_k\le\eps$, $0\le k\le n$.
\par If $a\le\phi^*+\psi^*<b$ and $a_k\le\phi^*<a_{k+1}$, then
$a-a_{k+1}<\psi^*<b-a_k$. Hence the set $G(a,b)=\{a\le\phi^*<b\}$
lies in the set
\[X\cup\left(\bigcup_{k=0}^{n-1}\left(\{a_k\le\phi^*<a_{k+1}\}\cap
\{a-a_{k+1}<\psi^*<b-a_k\}\right)\right).\]
\par Let $O'_k$ be open sets containing $\{a_k\le\phi^*<a_{k+1}\}$
and such that
\[\mu_j(\{\phi_j<a_k-\eps\}\cap O'_k)+\mu_j(\{\phi_j>a_{k+1}+\eps\}\cap
O'_k)<\frac\dl{n},\] while $O''_k$ be open sets containing
$\{a-a_{k+1}<\psi^*<b-a_k\}$ and such that
\[\mu_j(\{\psi_j<a-a_{k+1}-\eps\}\cap O''_k)+\mu_j(\{\psi_j>b-a_k+\eps\}\cap
O''_k)<\frac\dl{n}\] when $j$ is large. Let $O'''$ be an open set
containing $X$ and such that $\mu_j(O''')<\dl$ for large $j$. We
let $O_k=O'_k\cup O''_k$ and $O=O'''\cup\cup_{k=0}^{n-1}O_k$.
Clearly, $G(a,b)\sbs O$.
\par If $x\in O_k\cap K_j$,
$a_k-\eps\le\phi_j(x)<a_{k+1}+\eps$ and
$a-a_{k+1}-\eps\le\psi_j(x)\le b-a_k+\eps$, then
$a-2\eps<\phi_j(x)+\psi_j(x)<b+2\eps$. Thus the measure $\mu_j$
of those $x\in O_k\cap K_j$ for which
$\phi_j(x)+\psi_j(x)<a-2\eps$ or $\phi_j(x)+\psi_j(x)>b+2\eps$
does not exceed $2\dl/n$. Consequently,
\[\mu_j(\{\phi_j+\psi_j<a-2\eps\}\cap
O)+\mu_j(\{\phi_j+\psi_j>b+2\eps\}\cap O)<3\dl.\]
\end{proof}
\par The next theorem provides a convenient criterion for the existence of
limit values.
\bT\label{T:ecf} If $\phi$ has the strong limit values $\phi^*$ on $K$
then for every $\eps,\dl>0$ there is a function $f\in C(K)$ such that
$\mu(\{|f-\phi^*|>\eps\})<\dl$ and $\mu_j(\{|f-\phi_j|>\eps\})<\dl$ for
large $j$.
\par Moreover, if $\|\phi\|_{L^p(M)}<\infty$, $p\ge1$, then for every $\eps>0$
the function $f$ can be chosen so that
$\|f\|_{L^p(\mu_j)}<\|\phi\|_{L^p(M)}+\eps$ for large $j$ if $p<\infty$
and $\|f\|\le\|\phi\|_{L^\infty(M)}+\eps$ if $p=\infty$.
\par If, additionally, the measures $\{\mu_j\}$ converge
weak-$*$ to $\mu$ regularly and for a sequence $\phi\in L^p(M)$ and a
$\mu$-measurable function $\psi$ on $K_0$ and for every $\eps,\dl>0$
there is a sequence $f$ with strong limit values $f^*$ such that
$\mu(\{|f^*-\psi|>\eps\})<\dl$ and $\mu_j(\{|f_j-\phi_j|>\eps\})<\dl$
for large $j$, then $\phi$ has the strong limit values $\phi^*=\psi$ on
$K$.
\eT
\begin{proof} Suppose that $\phi$ has the strong limit values
$\phi^*$ on $K$. For any $\dl>0$ we can find $b>0$ such that
$\mu(\{\phi^*>b\})+\mu(\{\phi^*<-b\})<\dl$. Let
$X=\{x\in\supp\mu:\,-b\le\phi^*(x)\le b\}$. We fix $\eps>0$ and
cover the interval $[-b,b]$ with intervals with end points
$a_k=-b+k\eps$, $0\le k\le [2b/\eps]+1=n$ and for each $k$ take
a closed set $X_k\sbs G(a_{k-1},a_k)$ such that
$\sum\mu(X_k)>\mu(X)-\dl$. Then we select disjoint open sets
$O'_k$ containing $X_k$ such that (\ref{e:bv}) holds for the
given $\eps$ and $c_k\dl/n$, where $c_k=\min\{1,|a_k|^{-p}\}$.
\par The next step is to pick up continuous functions
$f_k$ on $K$ equal to 0 on $K\sm O'_k$, equal to $a_k$ on $X_k$
and taking values between 0 and $a_k$ elsewhere. Let $f=\sum
f_k$, $O_k=\{|f_k-a_k|<\eps\}$ and $O=\cup O_k$. Clearly,
$\mu(\{|f-\phi^*|>2\eps\})<3\dl$.
\par Note that $\mu_j(K)<\mu(K)+\dl$ for large $j$. Also for large $j$
we have
\[\mu_j(O)=\mu_j(\cup O_k)>\mu(\cup X_k)-\dl>\mu(X)-2\dl>\mu(K)-3\dl>
\mu_j(K)-4\dl.\]
\par Let $A_{jk}=\{x\in K_j\cap
O_k:\,|f(x)-\phi_j(x)|>3\eps\}$. Since
\[|f-\phi_j|\le|f-a_k|+|\phi_j-a_k|\] and $|f-a_k|<\eps$, we see that
if $x\in A_{jk}$, then $|\phi_j(x)-a_k|>2\eps$, i.e., either
$\phi_j(x)<a_k-2\eps$ or $\phi_j(x)>a_k+2\eps=a_{k+1}+\eps$. By
the choice of the sets $O_k$, the measure
$\mu_j(A_{jk})<c_k\dl/n$. Hence,
$\mu_j(\{|f-\phi_j|>3\eps\}\cap O\}<\dl$. But
$\mu_j(O)>\mu_j(K)-4\dl$. Hence $\mu_j(K\sm O)<3\dl$ and this
shows that $\mu_j(\{|f-\phi_j|>3\eps\}\}<4\dl$.
\par If $\|\phi\|_{L^\infty(M)}=A<\infty$, then the very
construction of $f$ shows that uniform norm of $f$ does not exceed
$A+\eps$. If $\|\phi\|_{L^p(M)}=A<\infty$, $1\le p<\infty$, then we
note that
\[\il K|f|^p\,d\mu_j=\sum\il{O'_k}|f|^p\,d\mu_j\]
and $|f|\le|a_k|$ on $O'_k$ while the measure of those $x\in
O'_k\cap K_j$, where $\phi_j(x)<a_k-\eps$ or
$\phi_j(x)>a_k+2\eps$, does not exceed $c_k\dl/n$ when $j$ is
large. Hence the measure of the set $B_{kj}=\{x\in O'_k\cap
K_j:\,|f(x)|>|\phi_j(x)|+2\eps\}$, does not exceed $c_k\dl/n$.
Thus
\[\il{O'_k}|f|^p\,d\mu_j\le\il{O'_k}(|\phi_j|+2\eps)^p\,d\mu_j+|a_k|^pc_k\dl/n
\le \il{O'_k}(|\phi_j|^p+2\eps)\,d\mu_j+\dl/n.\] Hence,
\[\il K|f|^p\,d\mu_j\le\il K(|\phi_j|+2\eps)^p\,d\mu_j+\dl.\]
By Minkowski's inequality
\[\il K(|\phi_j|+2\eps)^p\,d\mu_j\le
\left(\|\phi_j\|_{L^p(K_j,\mu_j)}+2\eps\mu_j^{1/p}(K)\right)^p.\]
Taking $\eps$ and $\dl$ sufficiently small we get the proof of the
statement.
\par To prove the last statement for given $\eps,\dl>0$ we
find a sequence $f$ with strong limit values $f^*$ satisfying the
conditions of the theorem. If $G(a,b)=\{a\le\psi<b\}$ and
$X=G(a,b)\sm\{a-\eps<f^*<b+\eps\}$, then $\mu(X)<\dl$. By the
regularity of the weak-$*$ convergence we can take an open set $O'$
such that $X\sbs O'$ and $\mu_j(O')<\dl$ for large $j$.
\par Let us take an open set $O''$ containing $\{a-\eps<f^*<b+\eps\}$
such that
\[\mu_j(\{f_j<a-2\eps\}\cap O'')+\mu_j(\{f_j>b+2\eps\}\cap
O'')<\dl.\] If $O=O''\cup O'$, then $G(a,b)\sbs O$ and
it is easy to see that
\[\mu_j(\{\phi_j<a-3\eps\}\cap O)+\mu_j(\{\phi_j>b+3\eps\}\cap
O)<2\dl.\]
\end{proof}
\par It is not true that the weak limit values are equal to the strong limit
values even when both do exist. For example, let all $\mu_j$ be
equal to the Lebesgue measure $\mu$ on $[0,1]$. Surround the
points $k/j$, $0\le k\le j$, by intervals of equal size and of
total length $1/j$. Let $\phi_j$ be equal to $j$ on these
intervals and $0$ outside. Then $\phi^*\equiv0$, while
$\phi_*\equiv1$.
\par However, under additional integrability assumptions, the strong
limit values correlate in the right way with powers and products
of sequences even when one of the factors has only the weak limit
values as the following chain of results shows.
\bT\label{T:cpp} Suppose that a sequence $\phi\in\A(M)$ has the strong limit
values $\phi^*$.
\be
\item If $\|\phi\|_{L^p(M)}=A<\infty$, $p>1$, then
    $\|\phi^*\|_{L^p(K,\mu)}\le A$ and for any $s$, $1\le s<p$,
    the measures $|\phi|^s\mu_j$ converge weak-$*$ to
    $|\phi^*|^s\mu$.
\item If $\|\phi\|_{L^p(M)}<\infty$, $p>1$, and
    a sequence $\psi\in\A^s(M)$, where $s>q$,
    $1/p+1/q=1$, then the sequence $\{\phi_j\psi_j\mu_j\}$ converges
    weak-$*$ to $\phi^*\psi_*\mu$.
\ee
\eT
\begin{proof} 1) Let us take a sequence of $\eps_n,\dl_n\sea0$. For
$\dl_n$ and $\eps_n$ let $f_n$ be a function from Theorem \ref{T:ecf}
with $\|f_n\|_{L^p(M)}\le A+\dl_n$. If $h\in C(K)$ then
\[\il Khf_n\,d\mu=\lim_{j\to\infty}\il Khf_n\,d\mu_j\le (A+\dl_n)\|h\|_{L^q(K,\mu)}.\]
Hence $f_n\in L^p(K,\mu)$ and $\|f_n\|_{L^p(K,\mu)}\le A+\dl_n$.
\par  The non-negative functions $|f_n|^p$ converge in measure to
$|\phi^*|^p$ and $\|f_n\|_{L^p(K,\mu)}\le A+\dl_n$. Hence, due to the
uniform integrability, $\|\phi^*\|_{L^p(K,\mu)}\le A$.
\par Let $K_{nj}=\{x\in K_j:\,|f_n(x)-\phi_j(x)|<\eps_n\}$. Firstly, since
$\mu_j(K\sm K_{nj})<\dl_n$ for large $j$,
\[\limsup_{j\to\infty}\il{K\sm K_{nj}}|\phi_j|^s\,d\mu_j\le
A^{s/p}\dl_n^{1-s/p}\]  and
\[\lim_{j\to\infty}\il{K\sm
K_{nj}}|f_n|^s\,d\mu_j\le (A+\dl_n)^{s/p}\dl_n^{1-s/p}.\] Secondly, if
$|x-y|<\eps_n$ then $||x|^s-|y|^s|\le s(|y|+\eps_n)^{s-1}\eps_n$.
\par If a function $h\in C(K)$, $\|h\|=1$, then
\begin{equation}\begin{aligned}
&\left|\il Kh(|\phi_j|^s-|f_n|^s)\,d\mu_j\right|= \left|\il
{K_{nj}}h(|\phi_j|^s-|f_n|^s)\,d\mu_j+\il {K\sm
K_{nj}}h(|\phi_j|^s-|f_n|^s)\,d\mu_j\right|\notag\\&\le
s\|\phi+\eps_n\|^{s-1}_{L^s(K_j,\mu_j)}\mu^{1/s}_j(K)\eps_n+2(A+\dl_n)^{s/p}\dl_n^{1-s/p}=
C_n,
\end{aligned}\end{equation}
where $C_n\to0$ as $n\to\infty$. Since
\[\lim_{j\to\infty}\il Kh|f_n|^s\,d\mu_j=\il Kh|f_n|^s\,d\mu,\] we see
that
\begin{equation}\begin{aligned}
\il Kh|f_n|^s\,d\mu-2C_n&\le \liminf_{j\to\infty}\il
Kh|\phi_j|^s\,d\mu_j\notag\\&\le \limsup_{j\to\infty}\il
Kh|\phi_j|^s\,d\mu_j\le\il Kh|f_n|^s\,d\mu+2C_n
\end{aligned}\end{equation}
for large $n$.
\par Let $K_n=\{x\in K_0:\,|f_n(x)-\phi^*(x)|<\eps_n\}$. Then estimates
similar to the used above give that
\[\left|\il Kh(|\phi^*|^s-|f_n|^s)\,d\mu\right|\le B_n,\]
where $B_n\to0$ as $n\to\infty$. This shows that
\[\lim_{n\to\infty}\il Kh|f_n|^s\,d\mu=\il Kh|\phi^*|^s\,d\mu.\]
Hence
\[\lim_{j\to\infty}\il Kh|\phi_j|^s\,d\mu_j\] exists and is
equal to
\[\il Kh|\phi^*|^s\,d\mu.\]
Thus $(|\phi|^s)_*=|\phi^*|^s$.
\par 2) The proof proceeds at the same style and with the same
notation. As before
\begin{equation}\begin{aligned}
&\left|\il Kh\psi_j(\phi_j-f_n)\,d\mu_j\right|=\left|\il
{K_{nj}}h\psi_j(\phi_j-f_n)\,d\mu_j+\il {K\sm
K_{nj}}h\psi_j(\phi_j-f_n)\,d\mu_j\right|\notag\\&\le
\eps_n\|\psi_j\|_{L^1(K,\mu_j})+\|\phi_j-f_n\|_{L^p(K,\mu_j)}\left(\il{K\sm
K_{nj}}|\psi_j|^q\,d\mu_j\right)^{1/q}.
\end{aligned}\end{equation}
\par Note that
\[\left(\il{K\sm
K_{nj}}|\psi_j|^q\,d\mu_j\right)^{1/q}\le\|\psi_j\|_{L^s(K,\mu_j)}\dl^{1/q-1/s}_n.\]
Since
\[\lim_{j\to\infty}\il Kh\psi_jf_n\,d\mu_j=\il Khf_n\psi_*\,d\mu,\] we see
that
\begin{equation}\begin{aligned}
\il Khf_n\psi_*\,d\mu-C_n&\le \liminf_{j\to\infty}\il
Kh\phi_j\psi_j\,d\mu_j\notag\\&\le \limsup_{j\to\infty}\il
Kh\phi_j\psi_j\,d\mu_j\le\il Khf_n\psi_*\,d\mu+C_n
\end{aligned}\end{equation}
for large $n$ and $C_n\to0$ as $n\to\infty$.
\par Estimates similar to the used above give that
\[\left|\il Kh(\phi^*\psi_*-f_n\psi_*)\,d\mu\right|\le B_n,\]
where $B_n\to0$ as $n\to\infty$. This shows that
\[\lim_{n\to\infty}\il Khf_n\psi_*\,d\mu=\il Kh\phi^*\psi_*\,d\mu.\]
Hence
\[\lim_{j\to\infty}\il Kh\phi_j\psi_j\,d\mu_j\] exists and is
equal to
\[\il Kh\phi^*\psi_*\,d\mu.\]
\end{proof}
\par Now we came to the main theorem of this section. It claims
that a sequence has the strong limit values when the inequality
proved in Theorem \ref{T:co} turns into an equality.
\bT\label{T:esbf} Let $\|\phi\|_{L^p(M)}\in\A^p(M)$ for some $p>1$, and
the measures $\{|\phi_j|^p\mu_j\}$ converge weak-$*$ to $\nu$.
If
\[\nu(K)=\il K|\phi_*|^p\,d\mu\]
then the sequence $\phi$ has the strong limit values equal to $\phi_*$.
\eT
\begin{proof} By Theorem \ref{T:co} $\nu=|\phi_*|^p\,d\mu$.
\par Let $f^+=\max\{f,0\}$ and $f^-=\max\{-f,0\}$. Since the
measures $\phi_j^+\mu_j$, $\phi_j^+\mu_j$, $(\phi_j^+)^p\mu_j$ and
$(\phi_j^-)^p\mu_j$ are uniformly bounded, we can find a
subsequence $j_k$ such that all these measures converge weak-$*$.
\par If the measures $\phi_{j_k}^+\mu_{j_k}$ converge weak-$*$ to $\psi\mu$,
then $\psi\ge0$ and the measures $\phi_{j_k}^-\mu_{j_k}$
converge weak-$*$ to $(\psi-\phi_*)\mu$. Hence
$\psi-\phi_*\ge0$.
\par If the measures $(\phi_{j_k}^+)^p\mu_{j_k}$ converge weak-$*$ to $\nu_1$
and the measures $(\phi_{j_k}^-)^p\mu_{j_k}$ converge weak-$*$
to $\nu_2$, then $\nu_1+\nu_2=\nu$ and by Theorem \ref{T:co}
$\nu_1\ge\psi^p\mu$ and $\nu_2\ge(\psi-\phi_*)^p\mu$. Hence
\[|\phi_*|^p\mu=\nu\ge\psi^p\mu+(\psi-\phi_*)^p\mu.\]
But the inequality $|x|^p\ge y^p+(y-x)^p$, where $y\ge0$ and
$y-x\ge0$, has only solution $y=x$ when $x\ge0$ or $y=0$ when
$x<0$. Hence $\psi=(\phi_*)^+$ and $\psi-\phi^*=(\phi_*)^-$.
\par It follows that the measures $\phi_j^+\mu_j$, $\phi_j^-\mu_j$,
$(\phi_j^+)^p\mu_j$ and $(\phi_j^-)^p\mu_j$ converge weak-$*$
to $(\phi_*)^+\mu$, $(\phi_*)^-\mu$, $(\phi_*^+)^p\mu$ and
$(\phi_*^-)^p\mu$ respectively. Therefore, it is sufficient to
prove our theorem only when all $\phi_j\ge0$.
\par Let us fix $0<\eps<1$ and $0<\dl<\eps/4$. We take $b\ge a>0$
such that $b/a<1+\dl$ and the set
$X=G(a,b)=\{x\in\supp\mu:\,a\le\phi_*(x)<b\}$. Suppose that
$\mu(X)>0$. Then we take an open set $O_1\sbs K$ and a closed set
$C\sps K$ such that $C\sbs X\sbs O_1$, $\mu(C)>0$ and $\mu(O_1\sm
C)<\dl\mu(X)$. We find a non-negative continuous function $f$ on
$K$ equal to 0 on $K\sm O_1$, equal to 1 on $C$, taking values
between 0 and 1 elsewhere and such that $\mu(f)<(1+\dl)\mu(C)$,
\[\dashint_ Kf\phi_*^p\,d\mu<b^p(1+\dl)\]
and
\[b(1+\dl)>\dashint_ Kf\phi_*\,d\mu>a(1+\dl)^{-1}.\]
\par Let $\mu'=f\mu$ and $\mu_j'=f\mu_j$. Since
$\mu'_j$, $\phi_j\mu'_j$ and $\phi_j^p\mu'_j$ converge weak-$*$
to $\mu'$, $\phi_*\mu'$ and $\phi^p_*\mu'$ respectively, there
is $j_0$ such that $\mu'_j(K)<(1+\dl)\mu(C)$,
\[\dashint_ K\phi_j^p\,d\mu'_j<b^p(1+\dl)\]
and
\[b(1+\dl)>\dashint_ K\phi_j\,d\mu'_j=\bar\phi_j>a(1+\dl)^{-1}\]
when $j\ge j_0$. Hence
\[\left(\dashint_ K\phi_j\,d\mu'_j\right)^{-p}\dashint_
K\phi_j^p\,d\mu'_j<(1+\dl)^{1+p}b^pa^{-p}\le(1+\dl)^{1+2p}\le1+r\dl,\]
where $r$ depends only on $p$.
\par If $\phi_j\ge (1+\eps)b$ then
\[\phi_j/\bar\phi_j\ge(1+\eps)(1+\dl)^{-1}>1+\eps/2\]
since $\dl<\eps/4$ and $0<\eps<1$. If $\phi_j\le(1-\eps)a$ then
$\phi_j/\bar\phi_j\le 1-\eps/2$. By Lemma \ref{L:chi}
\begin{equation}\begin{aligned}
\mu_j'(\{\phi_j\ge
&(1+\eps)b\})+\mu_j'(\{\phi_j\le(1-\eps)a\})\notag\\&\le
\mu_j'(\{|\phi_j/\bar\phi_j-1|>\eps/2\})\le\al(p,\eps/2)r\dl\mu'_j(K).
\end{aligned}\end{equation}
\par Let $O_2=\{f>1-\dl\}$. Then $O_2$ is open, lies in $O_1$ and contains $C$.
Moreover,
\begin{equation}\begin{aligned}
&\mu_j(O_2\cap\{\phi_j\ge(1+\eps)b\})+\mu_j(O_2\cap\{\phi_j\le(1-\eps)a\})
\notag\\&\le
(1-\dl)^{-1}(\mu_j'(\{\phi_j\ge(1+\eps)b\})+\mu_j'(\{\phi_j\le(1-\eps)a\}))
\notag\\&\le
2\al(p,\eps/2)r\dl\mu'_j(K)\le4\al(p,\eps/2)r\dl\mu_j(X).
\end{aligned}\end{equation}
\par Finally we use the regularity of the weak-$*$ convergence to
take an open set $O_3$ containing $X\sm C$ such that
$\mu_j(O_3)<\dl\mu(X)$ when $j$ is sufficiently large. Let
$O=O_2\cup O_3$. Then
\[\mu_j(\{\phi_j<(1-\eps)a\}\cap O)+\mu_j(\{\phi_j>(1+\eps)b\}\cap O)<
(4\al(p,\eps/2)r+1)\dl\mu(X).\]
\par Thus we have shown that for every $0<\eps<1$, $0<\dl<\eps/4$,
$b\ge a>0$ and the set $X=G(a,b)$, $\mu(X)>0$, there are a
constant $c(p,\eps)>0$ and an open set $O$ containing $X$ such
that
\begin{equation}\label{e:i1}
\mu_j(\{\phi_j<(1-\eps)a\}\cap O)+\mu_j(\{\phi_j>(1+\eps)b\}\cap
O)<c(p,\eps)\dl\mu(X)\end{equation} provided $b/a<1+\dl$.
\par Now for every $\eps,\dl>0$, $0<\eps<1$, $0<\dl<\eps/4$, we will show
the existence of open sets satisfying (\ref{e:bv}) when $0<a\le
b<1$. For this we, firstly, note that if $\mu(X)=0$, then the
existence of an open set for which (\ref{e:bv}) holds for all
$\eps$ and $\dl$ follows from the regularity of the weak-$*$
convergence. If $a>0$ and $\mu(X)>0$, then we define $a_0=a$,
$a_k=(1+\dl/2)a_{k-1}$ while $a_k<b$ and the last $a_n=b$. Let
$X_k=X\cap G(a_{k-1},a_k)$. If $\mu(X_k)=0$ then we cover it by
an open set $O_k$ such that $\mu_j(O_k)<\dl/n$ and if
$\mu(X_k)>0$ then we cover it by an open set $O_k$ such that
\[\mu_j(\{\phi_j<(1-\eps)a_{k-1}\}\cap O_k)+\mu_j(\{\phi_j>(1+\eps)a_k\}\cap
O_k)<c(p,\eps)\dl\mu(X_k).\]
\par Note that $a_{k-1}\ge a$ so $(1-\eps)a_{k-1}\ge
(1-\eps)a>a-\eps$ because $0<a<1$. By the same token $a_k\le b$
so $(1+\eps)a_k\le(1+\eps)b<b+\eps$ because $0<b<1$. Hence, if
$\phi_j(x)<a-\eps$ then $\phi_j(x)<(1-\eps)a_{k-1}$ and if
$\phi_j(x)>b+\eps$ then $\phi_j(x)>(1+\eps)a_k$. Therefore,
\[\mu_j(\{\phi_j<a-\eps\}\cap O_k)+\mu_j(\{\phi_j>b+\eps\}\cap
O_k)<c(p,\eps)\dl\mu(X_k).\] If $O=\cup O_k$ then
\[\mu_j(\{\phi_j<a-\eps\}\cap O)+\mu_j(\{\phi_j>b+\eps\}\cap
O)<c(p,\eps)\dl\mu(X)+\dl.\] Fixing $\eps$ and then for any
$\dl'>0$ picking up $\dl$ so that $0<\dl<\eps/4$ and
$c(p,\eps)\dl\mu(X)+\dl<\dl'$ we will get a set $O$ satisfying
(\ref{e:bv}).
\par If there are no an upper bound on $b$, then we take a
constant $\beta>b$ and consider the sequence
$\psi=\{\phi_j/\beta\}$. Evidently, $\psi_*=\phi_*/\beta$ and
the measures $\{|\psi_j|^p\mu_j\}$ converge weak-$*$ to
$\nu_p=\beta^{-p}\nu$. Hence $\nu_p=\psi_*^p\mu$.
\par Fix $\eps,\dl>0$ and take $\eps'<\eps/\beta$. By the previous result
there is an open set $O$ on $K$ containing
$X=\{a/\beta\le\psi_*<b/\beta\}$ such that
\[\mu_j(\{\psi_j<a/\beta-\eps'\}\cap O)+\mu_j(\{\psi_j>b/\beta+\eps'\}\cap
O)<\dl\] or
\[\mu_j(\{\phi_j<a-\eps'\beta\}\cap O)+\mu_j(\{\phi_j>b+\eps'\beta\}\cap
O)<\dl.\] But $\eps'\beta<\eps$ and we see that
\[\mu_j(\{\phi_j<a-\eps\}\cap O)+\mu_j(\{\phi_j>b+\eps\}\cap
O)<\dl.\]
\par If $a=0$ then we take $0<\beta<b$ such that
\[\frac{2\beta\mu(G(a,b))}{(\beta+\eps)(1-\beta)}<\dl.\]
If $Y=G(0,\beta)$ and $\mu(Y)>0$, then we take a closed set
$C\sbs Y$ such that the set $Y\sm C$ has so small measure that
there is an open set $O'$ containing $Y\sm C$ and of measure
less than $\dl$. Then we take $f\in C(K)$ equal to 1 on $C$ and
taking values between 0 and 1 elsewhere and such that
\[\il Kf\phi_*\,d\mu<2\beta\mu(Y).\] Let $O''=\{f>1-\beta\}$.
For large $j$ we have
\[2\beta\mu(Y)>\il
Kf\phi_j\,d\mu_j\ge(\beta+\eps)(1-\beta)\mu_j(\{\phi_j\ge\beta+\eps\}\cap
O'').\] Thus
\[\mu_j(\{\phi\ge\beta+\eps\}\cap O'')\le\frac{2\beta\mu(Y)}
{(\beta+\eps)(1-\beta)}<\dl.\] Taking $O=O'\cup O''$ we see
that it contains $Y$ and has all required properties. The set
$X\sm Y\sbs G(\beta,b)$ can be handled by the previous step.
\end{proof}
\section{Sequences satisfying integral inequalities}\label{S:ssii}
\par It became clear from the previous sections that to have strong
limit values the weak limit values need to exercise some control
over the values of functions. The typical form of this control is
\begin{equation}\label{e:tc}
\phi(z)\le\il{K_0}P(z,\zeta)\phi_*(\zeta)\,d\mu(\zeta),
\end{equation}
where $P$ is some kernel.
\par So let $M=\{\mu_j\}$ be a sequence of measures on
$K$ converging weak-$*$ to a finite measure $\mu$ on $K$ such that
$\supp\mu_j\cap\supp\mu_0=\emptyset$. Let $P(z,\zeta)$ be a
non-negative Borel function on $K\times K_0$. We require that for all
$j$ and for each fixed $z\in K_j$ the function $P(z,\zeta)$ is bounded
on $K_0$. Let $\A^p(M,P)$ be the set of sequences of Borel functions
$\phi_j$ defined on $K_j$ which have weak limit values $\phi_*$,
$\|\phi\|_{L^p(M)}<\infty$ and (\ref{e:tc}) holds for all
$z\in\cup_{j=1}^\infty K_j$.
\par It is reasonable to request that the class $\A^p(M,P)$
contains the constants and this is equivalent to request that
\begin{equation}\label{e:cf}
\il{K_0}P(z,\zeta)\,d\mu(\zeta)=1
\end{equation}
for all $z\in\cup_{j=1}^\infty K_j$.
\par We will need the following lemma.
\bL\label{L:cl} Let $\mu$ be a regular Borel measure on a compact space $X$.
If a function $f\in L^1(X,\mu)$, the Borel functions $p_j$ on
$X$ are uniformly bounded and converge weak-$*$ to $1$, then
\[\lim_{j\to\infty}\il Xfp_j\,d\mu=\il Xf\,d\mu.\]
\eL
\begin{proof} By Lusin's theorem for every $\eps>0$
we can find a closed set $Y\sbs X$ such that $f$ is continuous on
$Y$ and
\[\il{X\sm Y}|f|\,d\mu<\eps.\]
Let $B$ be the uniform norm of $f$ on $Y$. We take an open set
$O\sbs X$ such that $Y\sbs O$ and $\mu(O\sm Y)<\eps/B$. Then we
extend $f$ from $Y$ to $X$ as a continuous function $h$ such that
$\|h\|_X=B$ and $h=0$ on $X\sm O$. Note that
\[\il{X}|f-h|\,d\mu<3\eps.\]
\par Suppose that $|p_j|\le A<\infty$. Then
\[\left|\il X(f-h)p_j\,d\mu\right|\le\il
X|(f-h)p_j|\,d\mu\le3A\eps.\]Thus
\[\il Xhp_j\,d\mu-3A\eps\le\il Xfp_j\,d\mu\le\il Xhp_j\,d\mu+3A\eps.\]
Letting $j$ go to $\infty$ we get
\[\il Xh\,d\mu-3A\eps\le\liminf_{j\to\infty}\il Xfp_j\,d\mu\le
\limsup_{j\to\infty}\il Xfp_j\,d\mu\le\il Xh\,d\mu+3A\eps.\] Thus
\[\lim_{j\to\infty}\il Xfp_j\,d\mu=\il Xf\,d\mu.\]
\end{proof}
\par In the following theorem we introduce an important condition on the
kernel $p$ which plays a big role in the theory. It helps to prove a
theorem establishing the existence of strong boundary values for
sequences in $\A^p(M,P)$ and provides important estimates. For the
future, if $h$ is a function on $K_0$ we set the sequence $Ph$
consisting of functions
\[(Ph)_j(z)=\il{K}P(z,\zeta)h(\zeta)\,d\mu(\zeta),\qquad z\in K_j.\]
\bT\label{T:amp} Suppose that the kernel $P(z,\zeta)$ satisfies
(\ref{e:cf}), the functions \[p_j(\zeta)=\il
{K_j}P(z,\zeta)\,d\mu_j(z)\] are uniformly bounded and converge
weak-$*$ to $1$ on $K_0$. Then:\be\item any sequence
$\{\phi_j\}\in\A^p(M,P)$, $p>1$, has strong limit values equal to
$\phi_*$ and $\|\phi\|_{L^p(M)}=\|\phi_*\|_{L^p(K)}$;
\item if $h\in L^p(K_0,\mu)$ then $\|Ph\|_{L^p(M)}\le\|h\|_{L^p(K_0)}$;
\item if the sequence $Ph$ has strong limit values equal to $h$
    for a dense set of functions $h$ in $L^p(K_0,\mu)$ then the same holds for any
    function in $L^p(K_0,\mu)$, $p>1$.
\ee
\eT
\begin{proof} 1) Let $\phi\in\A^p(M,P)$. For $z\in K_j$ we define
\[\psi_j(z)=\il{K_0}P(z,\zeta)\phi_*(\zeta)\,d\mu(\zeta).\]
By (\ref{e:tc}) $\psi_j\ge\phi_j$ for all $j$. Moreover,
\[\il {K_j}\psi_j(z)\,d\mu_j(z)=
\il {K_0}\phi_*(\zeta)\left(\il {K_j}P(z,\zeta)\,d\mu_j(z)\right)
\,d\mu(\zeta)=\il{K_0}\phi_*(\zeta)p_j(\zeta)\,d\mu(\zeta)\] and
\[\il {K_j}|\psi_j(z)|\,d\mu_j(z)\le
\il {K_0}|\phi_*(\zeta)|\left(\il {K_j}P(z,\zeta)\,d\mu_j(z)\right)
\,d\mu(\zeta)=\il{K_0}|\phi_*(\zeta)|p_j(\zeta)\,d\mu(\zeta).\] Hence
the norms $\|\psi_j\|_{L^1(K_j)}$ are uniformly bounded and we can take
a subsequence $\{\psi_{j_k}\mu_{j_k}\}$ converging weak-$*$ to a
measure $\nu$. Then $\nu\ge\phi_*\mu$ but by Lemma \ref{L:cl}
\[\il {K_0}\,d\nu=\lim_{k\to\infty}\il
{K_{j_k}}\psi_{j_k}\,d\mu_{j_k}=\lim_{k\to\infty}\il{K_0}\phi_*p_{j_k}\,d\mu=
\il{K_0}\phi_*\,d\mu.\] Thus $\nu=\phi_*\mu$ and the sequence
$\{\psi_j\mu_j\}$ converges weak-$*$ to $\phi_*\mu$.
\par Therefore, the non-negative sequence $\{\psi_j-\phi_j\}$
has the zero weak-$*$ limit values. By Chebyshev's inequality
the sequence $\{\psi_j-\phi_j\}$ has the strong boundary values
equal to 0.
\par Now by the definition of $\psi$, H\"older's inequality and
(\ref{e:cf})
\[\il {K_j}|\psi_j|^p\,d\mu_j=
\il  {K_j}\left|\il {K_0}P(z,\zeta)\phi_*(\zeta)\,d\mu(\zeta)\right|^p
\,d\mu_j(z)\le\il{K_0}|\phi_*|^pp_j\,d\mu.\] Thus
$\|\phi\|_{L^p(M)}\le\|\phi_*\|_{L^p(K)}$. But by Lemma \ref{L:bv1}
$\|\phi\|_{L^p(M)}\ge\|\phi_*\|_{L^p(K)}$. Hence
$\|\phi\|_{L^p(M)}=\|\phi_*\|_{L^p(K)}$.
\par As before we derive that the sequence $\{|\psi_j|^p\mu_j\}$ is
bounded and if a subsequence $\{|\psi_{j_k}|^p\mu_{j_k}\}$ converges
weak-$*$ to a measure $\nu_p$, then
\[\il K\,d\nu_p\le\il{K}|\phi_*|^p\,d\mu.\]
By Theorem \ref{T:co} $\nu_p=|\phi_*|^p\mu$ and this implies that
the sequence $\{|\psi_j|^p\mu_{j_k}\}$ converges weak-$*$ to a
measure $\nu_p$. By Theorem \ref{T:esbf} the sequence $\{\psi_j\}$
has strong boundary values equal to $\phi_*$. Since
$\phi=(\phi-\psi)+\psi$ by Theorem \ref{T:pslf}(2) $\phi$ has the
strong limit values equal to $\phi_*$.
\par 2)
\begin{equation}\begin{aligned}
&\il{K_j}|(Ph)_j(z)|^p\,d\mu_j(z)=\il{K_j}
\left|\il{K_0}h(\zeta)P(z,\zeta)\,d\mu(\zeta)\right|^p\,d\mu_j(z)
\notag\\&\le
\il{K_j}\left(\il{K_0}|h(\zeta)|^pP(z,\zeta)\,d\mu(\zeta)
\left(\il{K_0}P(z,\zeta)\,d\mu(\zeta)\right)^{p/q}\right)\,d\mu_j(z)
\notag\\&=\il{K_j}\left(\il{K_0}|h(\zeta)|^pP(z,\zeta)\,d\mu(\zeta)
\right)\,d\mu_j(z)\notag\\&=\il{K_0}|h(\zeta)|^p
\left(\il{K_j}P(z,\zeta)\,d\mu(z)\right)\,d\mu(\zeta)=
\il{K_0}|h(\zeta)|^pp_j(\zeta)\,d\mu(\zeta).\notag
\end{aligned}\end{equation}
\par By Lemma \ref{L:cl} the last integrals converge to
$\|h\|^p_{L^p(K_0)}$.
\par 3) If $\phi\in L^p(K_0,\mu)$ then for every $\eps>0$ there is function $h\in
L^p(K_0,\mu)$ such that $Ph$ has strong limit values $h$ and
$\|h-\phi\|_{ L^p(K_0,\mu)}<\eps$. By 2)
$\|P(h-\phi)\|_{L^p(M)}\le\eps$. By the last statement in Theorem
\ref{T:ecf} the sequence $P\phi$ has strong limit values equal to
$\phi$.
\end{proof}
\par So we see the importance of the condition on the functions $p_j$
imposed in Theorem \ref{T:amp}. Let us list two important cases when
these conditions hold at least partially.
\bT\label{T:edf} Suppose that the kernel $P(z,\zeta)$ satisfies
(\ref{e:cf}).
\par 1)If the measures $P(z_j,\zeta)\mu(\zeta)$ converge
weak-$*$ to $\dl_{\zeta_0}$ when the sequence $\{z_j\}\sbs K$ converges
to $\zeta_0\in K_0$, then the functions $p_j(\zeta)$ converge weak-$*$
to $1$ on $K_0$. Moreover, if $h\in C(K_0)$ then the functions $(Ph)_j$
converge uniformly to $h$.
\par 2) If the measures $P(z,\zeta_j)\mu_j(z)$ converge weak-$*$
to $\dl_\zeta$ when the sequence $\{\zeta_j\}\sbs K_0$ converges to
$\zeta\in K_0$, then the functions $p_j(\zeta)$ are uniformly bounded
and converge weak-$*$ to $1$ on $K_0$. Moreover, if $h\in L^p(\mu)$,
$p>1$, then the sequence $Ph$ belongs to $\A^p(M,P)$ and has weak (and,
consequently, strong) limit values equal to $h$.
\eT
\begin{proof} 1) Let $h\in C(K_0)$. Let us show that for every $\eps>0$
there is $j_0$ and $\dl>0$ such that $|(Ph)_j(z)-h(\zeta)|<\eps$ when
$j\ge j_0$ and $|z-\zeta|<\dl$. Due to the compactness of $K$ the
negation of this statement means that there is $\zeta_0\in K_0$ and a
sequence $z_{j_k}\in K_{j_k}$ such that $|z_{j_k}-\zeta_0|\to0$ as
$k\to\infty$ but $|(Ph)_j(z_{j_k})-h(\zeta_0)|>\eps>0$. But it clearly
contradicts to the imposed condition on the measures
$P(z_j,\zeta)\mu(\zeta)$.
\par Hence, if $H$ is a continuous extension of $h$ to $K$ we may
assume that $|(Ph)_j(z)-H(z)|<\eps$ when $j\ge j_0$. Therefore,
\[\il{K_0}h(z)\,d\mu(\zeta)=\lim_{j\to\infty}\il{K_j}H(z)\,d\mu_j(z)=
\lim_{j\to\infty}\il{K_j}(Ph)_j(z)\,d\mu_j(z).\] But
\[\il{K_j}(Ph)_j(z)\,d\mu_j(z)=
\il{K_0}h(\zeta)\left(\il{K_j}P(z,\zeta)\,d\mu_j(z)\right)\,d\mu(\zeta)=
\il{K_0}h(\zeta)p_j(\zeta)\,d\mu(\zeta).\]
Hence,
\[\lim_{j\to\infty}\il{K_0}h(\zeta)p_j(\zeta)\,d\mu(\zeta)=\il{K_0}h(z)\,d\mu(\zeta)\]
and we are done.
\par 2) Clearly, the functions $p_j(\zeta)$ converge to 1
pointwise on $K_0$. To show that they are uniformly bounded we suppose
that there is a sequence $\{j_k\}$ and points $\zeta_k\in\supp\mu$ such
that $p_{j_k}(\zeta_k)\to\infty$ as $k\to\infty$. Without any loss of
generality we may assume that points $\zeta_k$ converge to $\zeta_0$.
Since the measures $P(z,\zeta_k)\mu_j(z)$ converge weak-$*$ to
$\dl_{\zeta_0}$ we got a contradiction.
\par Now if $f\in C(K)$ then
\[\il{K_j}f(z)(Ph)_j(z)\,d\mu_j(z)=\il
Kh(\zeta)\left(\il{K_j}f(z)P(z,\zeta)\,d\mu_j(z)\right)\,d\mu(\zeta).\]
The functions
\[F_j(\zeta)=\il{K_j}f(z)P(z,\zeta)\,d\mu_j(z)\] are uniformly
bounded and converge to $f(\zeta)$ pointwise. Hence
\[\lim_{j\to\infty}\il{K_j}f(z)(Ph)_j(z)\,d\mu_j(z)=\il
Kh(\zeta)f(\zeta)\,d\mu(\zeta).\] Thus the weak limit values of $\phi$
is $h$.
\par Now
\[\il{K_j}|(Ph)_j|^p\,d\mu_j\le\il K|h|^pp_j\,d\mu\] and we see
that $\|Ph\|_{L^p(M)}\le\|h\|_{L^p(\mu)}$. Thus $Ph\in\A^p(M,P)$ and
has the strong limit values equal to $h$.
\end{proof}
\par Generally, the inequality (\ref{e:tc}) comes from similar
inequalities obtained on the interior of a domain. This process is
described as follows:
\bT\label{T:lpk} Suppose that for each $j$ there are Borel functions
$P_j(z,\zeta)$
defined on $\cup_{m=1}^{j-1} K_m\times K_j$ such that for each
$z\in K_m$ the functions $P_j(z,\cdot)$, $j>m$, are uniformly
bounded and have the strong limit values $P(z,\cdot)$ with respect
to $M=\{\mu_j\}$.
\par If $\phi\in\A^p(M)$, $p>1$, and
\[\phi_m(z)\le\il{K_j}P_j(z,\zeta)\phi_j(\zeta)\,d\mu_j(\zeta)\]
for each $z\in\cup_{m=1}^{j-1}K_m$, then $\phi\in\A^p(M,P)$.
\eT
\begin{proof} Let us fix some $z\in K_m$ and let
$q_j(\zeta)=P_j(z,\zeta)$. The sequence $\{q_j\}\in L^\infty(M)$
and the sequence $\phi\in L^p(M)$, $p>1$. By Theorem
\ref{T:cpp}(2) the sequence $\{\phi_jq_j\mu_j\}$ converges
weak-$*$ to $\phi_*P(z,\cdot)\mu$. Hence
\[\phi(z)\le\il{K_0}P(z,\zeta)\phi_*(\zeta)\,d\mu(\zeta).\]
\end{proof}
\section{Boundary values}\label{S:bv}
\par These results can be applied to the theory of boundary values
in the following manner. Suppose that $D$ is a bounded domain in
$\mathbb R^n$ exhausted by domains $D_r$, $r<0$, such that $D_r\sbs\sbs
D_s$ when $s>r$ and $D=D_0=\cup_{r<0}D_r$. Suppose also that there are
measures $\mu_r$ supported by $S_r=\bd D_r$ converging weak-$*$ in
$C^*(\ovr D)$ as $r\to0^-$ to a finite measure $\mu$ supported by $\bd
D$. Let $K_0=\supp\mu$.
\par We say that a function $u$ on $D$ has {\it boundary values
with respect to measures $\mu_r$} if it has strong limit values with
respect to $M=\{\mu_{r_j}\}$ for any sequence $r_j\nea0$ and these
strong limit values do not depend on the choice of a sequence.
\par For $p\ge0$ and a continuous function $u$ on $D$ we define
\[\|u\|_p=\limsup_{r\to0^-}\il{S_r}|u|^p\,d\mu_r.\]
\par We assume that a Borel function $P(z,\zeta)$ is
defined on $D\times\bd D$.  We will require
the kernel $P(z,\zeta)$ to satisfy the following conditions:
\be
\item[(P1)] \[\il{\bd D}P(z,\zeta)\,d\mu(\zeta)=1;\]
\item[(P2)] for some $p>0$ and every $h\in L^p(K_0,\mu)$  the function
\[Ph(z)=\il{\bd D}h(z)P(z,\zeta)\,d\mu(\zeta)\]
has boundary values equal to $h$ and $\|Ph\|_p\le\|h\|_{L^p(K_0,\mu)}$.
\ee
\par To check the most complicated condition (P2) we note that by Theorem
\ref{T:amp} it holds if we assume that:\be
\item for every sequence $r_j\nea0$ the functions
\[p_{r_j}(\zeta)=\il{S_{r_j}}P(z,\zeta)\,d\mu_{r_j}(z)\] are uniformly
bounded and converge weak-$*$ to $1$ on $K_0$;
\item the set of $h\in L^p(K_0,\mu)$ such that the function $Ph$ has
    strong boundary values equal to $h$ with respect to any sequence
    $M=\{\mu_{r_j}\}$ is dense in $L^p(K_0,\mu)$.
\ee
The list of cases when (1) holds is given in Theorem \ref{T:edf}. The second
condition will follow if, for example, it holds for functions in $C(K_0)$.
\bT\label{T:bv} Let $D$ be a domain in $\mathbb R^n$. Let $\F^p$, $p>1$,
be a set of continuous functions $u$ on $D$ and let $P$ be a Borel function
on $D\times\bd D$ such that:
\be\item the kernel $P$ satisfies conditions (P1)--(P2);
\item $\|u\|_p<\infty$ for all $u\in\F^p$;
\item if $u\in\F^p$ and the functions $u|_{S_{r_j}}$ have weak limit
    values $\psi$ for some sequence $r_j\nea0$, then $u\le P\psi$ on $D$.
\ee
\par If $u\in\F^p$ then it has boundary values.
\eT
\begin{proof} Suppose that the functions $u_{r_j}$ have weak limit values
$\psi$ with respect to $M=\{\mu_{r_j}\}$ for some sequence $r_j\nea0$.
By Lemma \ref{L:bv1} $\psi\in L^p(K_0,\mu)$ and by (3) $u\le P\psi$ on
$D$. Hence the function $v=u-P\psi$ is non-positive and by (P2) has
weak boundary values with respect to $M$ equal to 0. Thus $v$ has
strong limit values with respect to $M$ equal to 0 and we see that $u$
has strong limit values with respect to $M=\{\mu_{r_j}\}$ equal to
$\psi$.
\par Now suppose that $u$ has strong
limit values $\phi$ with respect to $N=\{\mu_{t_j}\}$, $t_j\nea0$. The
functions $v_j=v|_{S_{t_j}}$ are non-positive and $\|v\|_p<\infty$.
Hence there is a subsequence $N'=\{\mu_{t_{j_k}}\}$ such that $v$ has
weak limit values $v_*$ with respect to $N'$. Clearly, $v_*\le0$. Thus
by (P2)
\[\phi=\{u|_{S_{t_{j_k}}}\}_*=\{(P\psi+v)|_{S_{t_{j_k}}}\}_*=\psi+v_*\] and we see
that $\phi\le\psi$. By the symmetry $\psi\le\phi$ and we see that
$\phi=\psi$.
\par If $u\in \F^p$ and $r_j\to0^-$ then by the first part of the proof all
subsequences of $\{u_{r_j}\}$ which have weak limit values have strong
limit values. By the second part of the proof these values don't depend
on the subsequence. Hence the sequence $\{u_{r_j}\}$ has weak limit
values and, consequently, strong limit values which do not depend on
the sequence. Hence, $u$ has boundary values.
\end{proof}
\par In this theorem the third condition is the most difficult to verify
(provided (P2) is checked). The only strategy for verification we know
is the following. Suppose that for each $r$ the kernels $P_r(z,\zeta)$
are defined on $D_r\times S_r$ and are Borel functions.  We will
require the kernels $P_r(z,\zeta)$ to satisfy the following conditions:
\be\item for each $z\in D$ the function $P_r(z,\zeta)$ is
    non-negative and uniformly bounded on $S_r$ when $r$ is close to 0;
\item $P_ru\ge u$ on $D_r$ for every $u\in\F^p$;
\item for each $z\in D$ the kernel $P(z,\cdot)$ is the strong limit
    values of the functions $P_{r_j}(z,\cdot)$ with respect to any sequence
    $M=\{\mu_{r_j}\}$, $r_j\nea0$.
\ee
\bT\label{T:v3} If the conditions above hold, $u\in\F^p$ and the functions
$u|_{S_{r_j}}$ have weak limit values $u_*$ for some sequence
$r_j\nea0$, then $u\le Pu_*$ on $D$ and, consequently, $u$ has boundary
values.
\eT
\begin{proof} The proof repeats the proof of Theorem
\ref{T:lpk}. Let us fix some $z\in D$ and let
$q_j(\zeta)=P_{r_j}(z,\zeta)$. Consider $M=\{\mu_{r_j}\}$. Then the sequence
$\{u_j=u|_{S_{r_j}}\}$ belongs to $L^p(M)$, $p>1$, and the sequence
$\{q_j\}\in L^\infty(M)$. By Theorem \ref{T:cpp}(2) the sequence
$\{u_jq_j\mu_{r_j}\}$ converges weak-$*$ to $u_*P(z,\cdot)\mu$. Hence
\[u(z)\le\il{\bd D}P(z,\zeta)u_*(\zeta)\,d\mu(\zeta).\]
\end{proof}
\par To check the third condition in the list above the following
result can be helpful. Assuming that $P=P_0$, $\mu=\mu_0$ and $\bd
D=S_0$ we introduce the function
\[Q_r(z,w)=\il{S_r}P_r(z,\zeta)P_r(w,\zeta)\,d\mu_r(\zeta).\]
\bT\label{T:Q} Suppose that:\be\item for each $z\in D$ the functions
$P_r(z,\zeta)$ are non-negative and uniformly bounded on $S_r$ when $r$
is close to 0;
\item for each $z\in D$ the kernel $P(z,\cdot)$ is the weak limit
    values of the functions $P_{r_j}(z,\cdot)$ with respect to any sequence
    $M=\{\mu_{r_j}\}$, $r_j\nea0$;
\item $Q_0(z,z)=\lim_{r\to0^-}Q_r(z,z)$ for every $z\in D$.
\ee
Then for each $z\in D$ the kernel $P(z,\cdot)$ is the strong limit
values of the functions $P_{r_j}(z,\cdot)$ with respect to any sequence
$M=\{\mu_{r_j}\}$, $r_j\nea0$.
\eT
\begin{proof} Suppose that the sequence $\{P^2_{r_j}(z,\cdot)\}$ has
weak  limit values equal to $\phi$   with respect to a sequence
$M=\{\mu_{r_j}\}$. Then
\[\il{K_0}\phi\,d\mu=
\lim_{j\to\infty}\il{S_r}P^2_{r_j}(z,\zeta)\,d\mu_{r_j}(\zeta)=
\lim_{j\to\infty}Q_{r_j}(z,z)=Q_0(z,z)=\il{K_0}P^2_0(z,\zeta)\,d\mu(\zeta).\]
By Theorem \ref{T:esbf} $P(z,\cdot)$ is the strong limit values of the
functions $P_{r_j}(z,\cdot)$ with respect to $M$.
\par Since the measures $P^2_r(z,\cdot)\mu_r$ are uniformly bounded
any sequence $\{\mu_{r_j}\}$, $r_j\nea0$, has a weak-$*$ converging
subsequence. But weak-$*$ limits of these sequences coincide and,
therefore, for each $z\in D$ the kernel $P(z,\cdot)$ is the strong
limit values of the functions $P_{r_j}(z,\cdot)$ with respect to any
sequence $M=\{\mu_{r_j}\}$, $r_j\nea0$.
\end{proof}
\par As the following theorem shows the existence of boundary values allows us
to expand the integral representation formulas from subdomains to the
whole domain. The proof follows immediately from the second part of
Theorem \ref{T:cpp}.
\bT\label{T:ir} Suppose that for each $z\in D$ the functions
$P_r(z,\zeta)$ are non-negative and uniformly bounded on $S_r$ when $r$
is close to 0 and the kernel $P(z,\cdot)$ is the weak limit values of
the functions $P_{r_j}(z,\cdot)$ with respect to some sequence
$M=\{\mu_{r_j}\}$, $r_j\nea0$. Suppose also that $\H$ is a set of
continuous functions $u$ on $D$ such that for any $u\in\H$:\be\item
\[u(z)=\il{S_r}u(\zeta)P_r(z,\zeta)\,d\mu_r(\zeta)+L_r(z,u),\]
for all $z\in D_r$ and $r<0$;\item $\|u\|_p<\infty$ for some
$p>1$;\item there are boundary values $u^*$.\ee If $L_r(z,u)\le0$ and
$\lim_{r\to0^-}L_r(z,u)=L(z,u)$ exists for all $z\in D$, then
\[u(z)=\il{S_0}u^*(\zeta)P(z,\zeta)\,d\mu(\zeta)+L(z,u).\]
\eT
\par In the last statement we show that when functions from some class
have radial limits then they have boundary values. By radial limits we
mean the following: a continuous mapping $v$ of $K_0\times[0,\eps)$,
$\eps>0$, into $\mathbb R^n$ is called {\it radial} if:\be\item
$v(\zeta,0)=\zeta$;
\item there is $r_0<0$ such that for every $\zeta\in K_0$ the set
    $\{t:\,v(\zeta,t)\in S_r\}\ne\emptyset$ when $r_0>r>0$, and the
    functions $t_r(\zeta)=\inf\{t:\,v(\zeta,t)\in S_r\}$ are
    continuous and converging uniformly to 0 as $r\nea0$;
\item there is a constant $c>0$ such that if a Borel set $E\sbs K_0$
    and $E_r=\{z\in S_r:\,z=v(\zeta,t_r(\zeta)),\zeta\in E\}$, then
    $\limsup_{r\to0^-}\mu_r(E_r)\le c\mu(E)$;
\item if $A_r=\{z\in
    S_r:\,z=v(\zeta,t_r(\zeta)),\zeta\in K_0\}$
    then $\limsup_{r\to0^-}\mu_r(S_r\sm A_r)=0$.
\ee
We say that a function $u$ on $D$ has radial limits $\wtl u$ $\mu$-a.e.
with respect to a radial mapping $v$  if there is a function $\wtl u$
on $K_0$ such that $\lim_{t\to0}u(v(\zeta,t_r(\zeta)))$ exists
$\mu$-a.e and is equal to $\wtl u(\zeta)$ $\mu$-a.e.
\bT\label{T:rv} If a Borel function $u$ on $D$ has radial limits $\mu$-a.e. with
respect to a radial mapping $v$, then $u$ has boundary values equal to
$\wtl u$.
\eT
\begin{proof} First of all, we note that by Lemma \ref{L:rc} for any sequence
$r_j\nea0$ the measures $\mu_{r_j}$ converge to $\mu$ regularly. So we
can use the last part of Theorem \ref{T:ecf}. Given a function $u$ on
$D$ with radial limits $\wtl u$ we fix $\eps,\dl>0$ and take a
continuous function $f$ on $K_0$ such that $f=\wtl u$ on a set $E$ with
$\mu(E)>\mu(K_0)-\dl$. Let us denote by the same letter $f$ the
continuous extension of $f$ to $D$.
\par Since the functions $t_r$ converge uniformly to 0 there are
$r_0<0$ and a set $F\sbs K_0$ such that $\mu(F)>\mu(K_0)-\dl$,
\[\sup_{0<r<r_0}|\wtl u(\zeta)-u(v(\zeta,t_r(\zeta)))|<\eps\] and
\[\sup_{0<r<r_0}|\wtl f\zeta)-f(v(\zeta,t_r(\zeta)))|<\eps\]
for all $\zeta\in F$.
\par Let $G_r$ be the set of points $z$ in $S_r$ such that
$z=v(\zeta,t_r(\zeta))$ for some $\zeta\in F\cap E$. By the property
(3) of radial mappings we can find negative $r_1\ge r_0$ such that
$\mu_r(S_r\sm A_r)<\dl$ when $r>r_1$. By the property (2) of radial
mappings we can find negative $r_2\ge r_1$ such that $\mu(K_0\sm (E\cap
F))_r)<2c\dl$ when $r>r_2$. Hence
\[\mu_r(\{|u(\zeta)-f(\zeta)|>2\eps\})\le (2c+1)\dl.\]
Hence $u$ has boundary values equal to $\wtl u$.
\end{proof}
\section{Boundary values of harmonic functions}\label{S:bvhf}
\par The simplest case is the case of harmonic functions on
regular domains, i.e., bounded domains such that any continuous
function $\phi$ on their boundaries has a harmonic extension $h_\phi$
to the domain continuous up to the boundary and this extension
coincides with $\phi$ on the boundary.
\par Let $D$ be such a domain. For any $x\in D$ the evaluation $h_\phi(x)$
defined for $\phi\in C(\bd D)$ is a continuous linear functional by the
maximum principle and, consequently, there are harmonic measures $\mu_x$ on
$\bd D$ such that
\[h_\phi(x)=\il{\bd D}\phi\,d\mu_x.\]
\par For any $\zeta\in\bd D$ let $S_{\zeta,r}$ be the intersection of
$\bd D$ with the closed ball centered at $\zeta$ and of radius $r>0$
and let $\chi_{\zeta,r}$ be the characteristic function of
$S_{\zeta,r}$. Since $\chi_{\zeta,r}$ is the limit of a decreasing
sequence of functions in $C(\bd D)$, the function
$F_r(x,\zeta)=\mu_x(S_{\zeta,r})$ is harmonic on $D$ and, clearly, for
any $\zeta\in\bd D$ extends continuously as 0 on $\bd D\sm S_{\zeta,r}$
and as 1 on the relative interior of $S_{\zeta,r}$ in $\bd D$. We
define $F_r(x,\zeta)$ on $\ovr D\times\bd D$ letting it to be equal to
1 when $x\in S_{\zeta,r}$. So $F_r$ is a Borel function on $\ovr
D\times\bd D$.
\par Let us fix $x_0\in D$ and define
\[P(x,\zeta,r)=\frac{F_r(x,\zeta)}{F_r(x_0,\zeta)}.\]
\par By Harnack's inequality for every compact set $K\sbs D$
there is a constant $C(K)>0$ such that $C^{-1}(K)F_r(x,\zeta)\le
F_r(x_0,\zeta)\le C(K)F_r(x,\zeta)$ when $x,x_0\in K$ and $\zeta\in\bd
D$. Hence the function
\[P(x,\zeta)=\limsup_{r\to0^-}P(x,\zeta,r)\]
is Borel on $D\times\bd D$, $C^{-1}(K)\le P(x,\zeta)\le C(K)$ when
$x,x_0\in K$ and $\zeta\in\bd D$ and, again by Harnack's inequality,
$P(x,\zeta)$ is continuous in $x$. Thus $P(x,\zeta)$ is subharmonic in
$x$. Moreover, $\mu_x(\zeta)=P(x,\zeta)\mu_{x_0}(\zeta)$ for any $x\in
D$. Hence
\[h_\phi(x)=\il{\bd D}\phi(\zeta)P(x,\zeta)\,d\mu_{x_0}.\]
\par Let $V$ be the volume element on $D$. We cover $D$ by
a countable family of closed balls $B(x_j,r_j)$ centered at $x_j$ and
of radius $r_j>0$ such that $B(x_j,r_j)\sbs D$. For every $j$ and
$\zeta\in\bd D$ let
\[\Psi_j(x,\zeta)=\frac1{V(B(x_j,r_j))}\il{B(x_j,r_j)}P(y,\zeta)\,dV.\]
Then
\[h_\phi(x_j)=\frac1{V(B(x_j,r_j))}\il{B(x_j,r_j)}h_\phi(y)\,dV=
\il{\bd D}\phi(\zeta)\Psi_j(x_j,\zeta)\,d\mu_{x_0}.\] Since it is true
for every $\phi\in C(\bd D)$ we see that the equality
$\Psi_j(x_j,\zeta)=P(x_j,\zeta)$ holds for all $\zeta\in\bd D$ except
of a set $E_j$ with $\mu_{x_0}(E_j)=0$. But the subharmonicity of
$P(x,\zeta)$ yields that $P(y,\zeta)$ is harmonic on $B(x_j,r_j)$ for
all $\zeta\in\bd D\sm E_j$. Hence the function $P(x,\zeta)$ is harmonic
in $x$ for $\mu_{x_0}$-almost all $\zeta$ and is the Poisson kernel on
$D$ centered at $x_0$.
\par Now suppose that $D$ is exhausted by domains $\{D_r\}$, $r<0$.
Let $\mu_{r,x}$ and $P_r(x,\zeta)$ be the harmonic measure and the
Poisson kernel centered at $x_0$ respectively on $D_r$. Note that
$P(x_0,\zeta)=P_r(x_0,\zeta)=1$. Hence by Harnack's inequality for
every compact set $F\sbs D$ there is a constant $C(F)>0$ such that
$C^{-1}(F)\le P_r(x,\zeta)\le C(F)$ when $x\in F$.
\par Since for any $\phi\in C(\ovr D)$ and any $x\in D$
\[\lim_{r\to0^-}\il{S_r}
\phi(\zeta)P_r(x,\zeta)\,d\mu_{r,x_0}(\zeta)= \lim_{r\to0^-}
\il{S_r}h_\phi(\zeta)P_r(x,\zeta)\,d\mu_{r,x_0}(\zeta)=h_\phi(x),\] we see
that the measures $P_r(x,\zeta)\,d\mu_{r,x_0}(\zeta)$ converge weak-$*$
to $P(x,\zeta)\,d\mu_{x_0}(\zeta)$. In particular, since
$P_r(x_0,\zeta)=P(x_0,\zeta)\equiv1$ the  measures $\mu_{r,x_0}(\zeta)$
converge weak-$*$ to $\mu_{x_0}(\zeta)$.
\par In the future $p\ge1$ and we let $\mu_{r,x_0}=\mu_r$ and $\mu_{x_0}=\mu$.
We define the space $\rS^p(D)$ as the space of all continuous
subharmonic functions $u$ on $D$ such that
\[\|u\|^p_{\rS^p(D)}=\limsup_{r\to0^-}\il{S_r}|u|^p\,d\mu_r<\infty.\]
\par  The following theorem shows that any function in $L^p(\bd D,\mu_{x_0})$ is the trace of a function in $\rS^p(D)$ and the Poisson integral serves as restoring operator.
\bT\label{T:bvhf} Let $D$ be a regular domain in $\mathbb R^n$ and $p>1$.
For every function $\phi\in L^p(\bd D,\mu_{x_0})$ there is a unique
harmonic function $h_\phi\in\rS^p(D)$ with boundary values equal to
$\phi$ and $\|h_\phi\|_{\rS^p(D)}\le\|\phi\|_{L^p(\bd D,\mu_{x_0})}$.
Moreover, $h_\phi=P\phi$.
\eT
\begin{proof} Note that
\[\il{\bd D}P(x,\zeta)\,d\mu(\zeta)\equiv1\]
so the condition (P1) from Section \ref{S:bv} holds. Since $P(x,\zeta)$
is harmonic for $\mu$-almost all $\zeta$ we see that
\[p_r(\zeta)=\il{\bd D}P_r(x,\zeta)\,d\mu_r(\zeta)=P_r(x_0,\zeta)=1.\]
\par Moreover, if $h\in C(\bd D)$ then the function $Ph$ has boundary values
equal to $h$ with respect to any sequence $M=\{\mu_{r_j}\}$.  Since
$C(\bd D)$ is dense in $L^p(\bd D,\mu_{x_0})$, by Theorem \ref{T:amp}
for every $\phi\in L^p(\bd D,\mu_{x_0})$  the function
\[h_\phi(x)=P\phi(x)=\il{\bd D}\phi(x)P(x,\zeta)\,d\mu(\zeta)\]
has boundary values equal to $h$ with respect to any sequence
$M=\{\mu_{r_j}\}$ and $\|P\phi\|_p\le\|\phi\|_{L^p(\bd D,\mu_{x_0})}$.
Thus the condition (P2) holds also.
\par If $u\in\rS^p(D)$ is another function with  boundary values $\phi$,
then we take a sequence $\{r_j\nea0\}$. The weak limit of functions
$P_{r_j}(x,\zeta)$ with respect to $\{\mu_{r_j}\}$ is equal to
$P(x,\zeta)$ for all $x\in D$. By Theorem \ref{T:cpp}(2)
\[u(x)=\lim_{j\to\infty}\il{\bd D_{r_j}}u(y)P_r(x,y)\,d\mu_{r_j}(y)=
\il{\bd D}\phi(\zeta)P(x,\zeta)\,d\mu(\zeta).\] Thus $u=P\phi$.
\end{proof}
\par By the Riesz Representation Theorem every subharmonic function $u$ on
$D$ can be represented on $D_r$  as
\[u(x)=\il{S_r}u(y)P_r(x,y)\,d\mu_r(y)+\il {D_r}G_{D_r}(x,y)\,\Delta u(y),\]
where  $G_{D_r}(x,y)$ is the Green kernel on $D_r$. For functions in
$\rS^p(D)$ this property is characteristic for functions with boundary
values. \bT A function $u\in\rS^p(D)$ has the representation
\begin{equation}\label{e:rrf}
u(x)=\il{S}\phi(y)P(x,y)\,d\mu(y)+\il {D}G_{D}(x,y)\,\Delta u(y).
\end{equation}
if and only if it has boundary values equal to $\phi$.
\eT
\begin{proof}
\par If $u\in\rS^p(D)$ has boundary values $\phi$ then (\ref{e:rrf}) holds by
Theorem \ref{T:ir}.
\par In view of Theorem \ref{T:bvhf} to prove the converse it suffices to
show that the potential \[v(x)=\il {D}G_{D}(x,y)\,\Delta u(y)\] has
zero boundary value. For this since $v\le0$ it suffices to show that
$v$ has weak limit values equal to 0 with respect to any sequence
$M=\{\mu_{r_j}\}$, $r_j\nea 0$. \par Let $C\ge f\ge0$ be a continuous
function on $\ovr D$. Suppose that $r<0$ and $x_0\in D_r$. Then
\[f_r(y)=\il{S_r}f(x)G_{D}(x,y)\,d\mu_r(x)\ge C\il{S_r}G_{D}(x,y)\,d\mu_r(x).\]
Since $G_{D}(x,y)$ is subharmonic in $x$ when $y\in D_r$ and harmonic in $x$
on $D_r$ when $y\in D\sm\ovr D_r$,
\[\il{S_r}G_{D}(x,y)\,d\mu_r(x)\ge G_{D}(x_0,y).\] Moreover, for any $r_0<0$
and any $\eps>0$  there is $r_1<0$ such that if $r_1>r>0$ and $y\in
D_{r_0}$, then $G_{D}(x,y)>-\eps$ when $x\in D\sm\ovr D_{r_1}$. Thus
$f_r(y)\ge-C\eps$ on $D_{r_0}$ and $f_r(y)\ge-CG_D(x_0,y)$ on the rest
of the domain. Hence
\[\il{S_r}f(x)v(x)\,d\mu_r(x)=\il Df_r(y)\,\Delta u(y)\to0\] as $r\to0^-$
and we see that $v$ has weak limit values equal to 0 with respect to
any sequence $M=\{\mu_{r_j}\}$, $r_j\nea 0$.
\end{proof}
\section{Spaces of plurisubharmonic functions}\label{S:sopf}
\par Let $D$ be a hyperconvex domain in $\mathbb C^n$, i.e. it has a continuous
plurisubharmonic exhausting function $u$ equal to 0 on $\bd D$. We
assume that such functions can take $-\infty$ as their value.
\par On such domains for each $w\in D$ there is a unique continuous
plurisubharmonic function $G(z,w)=g_w(z)$ on $\ovr D\times D$ equal to
0 on $\bd D\times D$, satisfying the equation $(dd^cg_w)^n\equiv0$ on
$D\sm\{w\}$ and $g_w(z)-\log|z-w|=O(1)$ as $z\to w$. The function $g_w$
is called the pluricomplex Green function on $D$ with the pole at $w$.
\par For $w\in D$ and $r<0$ we let $g_{w,r}=\max\{g_w,r\}$. Following
\cite{D} we let
\[\mu_{w,r}=\frac1{(2\pi)^n}((dd^cg_{w,r}^n-\chi_{D\sm B_{w,r}}(dd^cg_w)^n),\]
where the domain $B_{w,r}=\{z\in D:\,g_w<r\}$. The measure $\mu_{w,r}$
is nonnegative and supported by $S_{w,r}=\{z\in D:\,g_w=r\}$. As $r\to
0^-$ the measures $\mu_{w,r}$ converge weak-$*$ to a measure $\mu_w$.
\par In \cite[Chapter 5]{D} Demailly had proved the
following Lelong--Jensen formula.
\bT\label{T:pir} Let $D$ be a hyperconvex domain and $w_0\in D$. There is a
non-negative
Borel function $P_D(w,\zeta)$ on $D\times\bd D$ such that
$\mu_w(\zeta)=P_D(w,\zeta)\mu_{w_0}(\zeta)$ and if $u$ is a continuous
function on $\ovr{D}$ plurisubharmonic on $D$, then
\[u(w)=
\il{\bd D}u(\zeta)P_D(w,\zeta)\,d\mu_{w_0}(\zeta)+
\frac1{(2\pi)^n}\il Dg_wdd^cu\wedge(dd^cg_w)^{n-1}.\]
\eT
\par In particular, if $u\equiv 1$ we see that
\[\il{\bd D}P_D(w,\zeta)\,d\mu(\zeta)=1.\]

\par To expand the space $PS^c(D)$ of continuous function on $\ovr{D}$
plurisubharmonic
on $D$ the following construction is natural (see \cite{PS} for more
details). Let us fix a point $w_0\in D$ and define the space
$PS^p_{w_0}(D)$ as the space of plurisubharmonic functions $u$ on $D$
such that
\[\limsup_{r\to0^-}\il{S_{w_0,r}}|u|^p\,d\mu_{w_0,r}<\infty.\]
We will denote the latter limit as $\|u\|_{PS^p_{w_0}(D)}$.
\par The domains $B_r=B_{w_0,r}$ are also hyperconvex. For $r<0$ we denote by
$g_{w,r}(z)=g_{B_r}(z,w)$ the pluricomplex Green function on
$B_r=B_{w_0,r}$ and let $\wtl\mu_{w,r}$ be the measure
generated by $g_{w,r}$ on $S_r=S_{w_0,r}$.
\bL\label{L:ucpm} For any $w\in D$ the measures $\wtl\mu_{w,r}$ converge weak-$*$
to $\mu_w$ as $r\to0^-$, i.e., $\wtl\mu_{w,r}(u)$ converge to
$\mu_w(u)$ as $r\to0^-$ when $u\in C(\ovr D)$.
\eL
\begin{proof} Since any function in $C(\ovr D)$ can be uniformly
approximated by $C^2$-functions on $\mathbb C^n$ which are differences
of plurisubharmonic $C^2$-functions, it suffices to prove this theorem
when $u$ is plurisubharmonic and $C^2$ on $\mathbb C^n$.
\par By Theorem \ref{T:ir}
\[u(w)=
\il{S_r}u(\zeta)\,d\wtl\mu_{w,r}(\zeta)+ \frac1{(2\pi)^n}\il
{B_r}g_{w,r}dd^cu\wedge(dd^cg_{w,r})^{n-1}.\] Since the functions
$g_{w,r}$ decrease when $r$ increases and converge to $g_w$ we see that
the volume integrals in the latter formula converge to
\[\il Dg_wdd^cu\wedge(dd^cg_w)^{n-1}.\] Hence $\wtl\mu_{w,r}(u)$ converges to
$\mu_w(u)$ as $r\to0^-$.
\end{proof}
\par  Let
$P_r=P_{B_r}$. Our first step is the following lemma.
\bL\label{L:brgf} Let $D$ be a hyperconvex domain in $\mathbb C^n$, let
$K$ be a compact set in $D$ and $w_0\in K$. There is a constant $c>0$
depending only on $K$ and $D$ such that  $c\wtl\mu_{w_1,r}\le
\mu_{w_2,r}(z)\le c^{-1}\wtl\mu_{w_1,r}$ for any points $w_1,w_2\in K$
when $r$ is sufficiently close to 0.
\par Consequently, the functions $P_r(z,w)$ are uniformly bounded on
$S_r$ when $w\in K$ and $r$ is sufficiently close to 0.
\eL
\begin{proof} Let us take $r_0$ such $K\sbs B_{r_0}=B_{w_0,r_0}$. By the
continuity of the pluricomplex Green functions there is a constant
$a>0$ such that $-a\le g_w\le-a^{-1}$ on $S_{r_0}$ for any $w$ in $K$.
Hence there is a constant $b>0$ such that $bg_{w_1}\le g_{w_2}\le
b^{-1}g_{w_1}$ for any $w_1,w_2\in K$. By the maximality of
pluricomplex Green functions this inequality holds on $D\sm B_{r_0}$.
\par The functions $g_{w,r}$ are continuous, decreasing in $r$ and
converging to $g_w$ pointwise. Hence they converge to $g_w$ uniformly
on compacta as $r\to0^-$ and we can find $r_1$ between $r_0$ and 0 such
that $g_w\le g_{w,r}<2g_w$ on $S_{r_0}$ for all $w\in K$. Thus there is
a constant $c>0$ such that $cg_{w_1,r}\le g_{w_2,r}\le c^{-1}g_{w_1,r}$
on $S_{r_0}$ for any $w_1,w_2\in K$ when $r$ is sufficiently close to
0. By the maximality of pluricomplex Green functions this inequality
holds on $B_r\sm B_{r_0}$. By \cite[Theorem 3.8]{D}
$c\wtl\mu_{w_1,r}\le \wtl\mu_{w_2,r}(z)\le c^{-1}\wtl\mu_{w_1,r}$.
\par Since $\wtl\mu_{w,r}(z)=P_r(w,z)\wtl\mu_{w_0,r}$ we see that the
functions $P_r(z,w)$ are uniformly bounded on $S_r$ when $w\in K$.
\end{proof}
\par As the result of two lemmas above we have the following corollary.
\bC Let $D$ be a hyperconvex domain and $w_0\in D$. If a function
$u\in PS^p_{w_0}(D)$, $p>1$, has boundary values $u^*$, then
\begin{equation}\label{e:irp}u(w)= \il{\bd
D}u^*(\zeta)P_D(w,\zeta)\,d\mu_{w_0}(\zeta)+ \frac1{(2\pi)^n}\il
Dg_wdd^cu\wedge(dd^cg_w)^{n-1}.
\end{equation}
\eC
\begin{proof} For each $z\in D$ the functions $P_r(z,\zeta)$ are non-negative
and uniformly bounded on $S_r$ when $r$ is close to 0 and the kernel
$P(z,\cdot)$ is the weak limit values of the functions
$P_{r_j}(z,\cdot)$ with respect to any sequence $M=\{\mu_{w_0,r_j}\}$,
$r_j\nea0$
\par The equation (\ref{e:irp}) holds if $D$ is replaced by $D_r$, $r<0$, and
$u^*$ replaced by $u$. Since the surface integrals in (\ref{e:irp})
stay bounded and volume integrals are negative and decrease as
$r\to0^-$, by Theorem \ref{T:ir} we have (\ref{e:irp}) on $D$.
\end{proof}
\par
\bT\label{T:irfp} Let $D$ be a strongly pseudoconvex domain with $C^2$
boundary and $w_0\in D$. If a function $u\in PS^p_{w_0}(D)$, $p>1$,
then it has boundary values $u^*$.
\eT
\begin{proof} Let $\rho$ be a defining function of $D$. It was proved in
\cite[Theorem 6.1]{D} that there are positive constants $C_1$ and $C_2$
such that
\[C_1(dd^c\rho)^{n-1}\wedge d^c\rho\le\mu_{w_0,r}\le
C_2(dd^c\rho)^{n-1}\wedge d^c\rho\] on $S_r$. The form
$(dd^c\rho)^{n-1}\wedge d^c\rho$ is continuous and strictly positive
and, therefore, there are positive constants $A_1$ and $A_2$ such that
$A_1\lm_r\le\mu_{w_0,r}\le A_2\lm_r$ on $S_r$, where $\lm_r$ is the
surface area.
\par Let $n(\zeta)$ be the inward normal vector at $\zeta\in\bd D$. It
follows from above that $\zeta+tn(\zeta)$ is a radial mapping. Since
$u$ is a subharmonic function on $D$ which has a harmonic majorant, by
\cite{S} $u$ has radial limits $\mu_{w_0}$-a.e. with respect to $n$. By
Theorem \ref{T:rv} $u$ has boundary values $u^*$.
\end{proof}
\par As an example of a non-smooth domain we can offer only a polydisk.
If $D=\mathbb D^n$, $w=(w_1,\dots,w_n)$ and $z=(z_1,\dots,z_n)$, then
(see \cite{D})
\[g_w(z)=\max\limits_{1\le j\le n}\log\left|\frac{z_j-w_j}{1-\bar w_jz_j}\right|.\]
\bT\label{T:irfpd} Let $D=\mathbb D^n$ and $w_0\in D$. If a function
$u\in PS^p_{w_0}(D)$, $p>1$, then it has boundary values $u^*$.
\eT
\begin{proof} We will prove it for $w_0=0$. Other cases are easily obtained
by biholomorphic transformations moving $w_0$ to 0. Since (see
\cite{D}) $\mu_r$ is supported by the Shilov boundary of $\mathbb
D_r^n$ and is equal to $(2\pi)^{-n}d\th_1\dots d\th_n$,
$z_j=re^{i\th_j}$, we see that if $u\in PS^p_0(D)$ and
$\th=(\th_1,\dots,\th_n)$, then $u\in PS^p_0(\mathbb D_\th)$, where
$\mathbb D_\th=\{z=\xi(e^{i\th_1},\dots,e^{i\th_n}),\xi\in\mathbb D\}$
for $\mu_0$-almost all $\th$.
\par If the radial mapping is $v(\zeta,t)=(1-t)\zeta$ then $u$ has radial
limits $\mu_0$-almost all $\th$ and by Theorem \ref{T:rv} $u$ has
boundary values $u^*$.
\end{proof}

\end{document}